\documentclass{amsart}
\usepackage{amssymb,latexsym}
\begin{document}

\title[Automorphism Groups of Groups of Order $16p^2$]{The Automorphism Groups\\ of the Groups of Orders $16p$ and $16p^2$}

\author{Elaine W. Becker}
\address{American Mathematical Society, Providence, Rhode Island}
\author{Walter Becker}
\address{266 Brian Drive, Warwick, Rhode Island 02886}
\maketitle

\begin{abstract}

            Results of the computation of the automorphism groups for
            the groups of orders $16p$ and $16p^{2}$ are given. In some cases
            it has not been possible to give as complete a set of results
            as was done previously for the case of groups of order $8p^2$.
            Problems arise for those groups of the form ($C_{p} \times  C_{p}$) @ $G$[16]
            that occur in the orders $p\equiv 1$ mod(8) and $p\equiv  7$ mod(8),
            where $G$[16] means any group of order 16.
            The groups $G$[16] in question are $C_{16}$, $D_{8}$, $QD_{8}$, and
            $Q_{4}$. For the other cases, explicit presentations are presented for
            the automorphism groups of the groups of orders 16$p$ and 16$p^2$.
\end{abstract}

\section {\underline{Introduction}}

            This is the second in a series of reports giving the results
      of a systematic computer study of the automorphism groups of finite
      groups of low order. Some background on the origin and extent of the
      work done to date can be found in reference \cite{1}. The first article
      in this series was devoted to the groups of order $8p^2$; see \cite{2}.

            The groups of order 16 and $16p$ for $p$ an odd prime have
      been known for almost a hundred years. See, for example, Lunn and Senior
      \cite{11} and the references therein. The automorphism groups
      of the groups of order $16p$ however have apparently not been computed.
      This report lists the groups of orders 16, $16p$, and $16p^2$ along with
      their automorphism groups. The automorphism groups for the groups of
      order 16 appear to have been known for a long time \cite{3}. The
      automorphism groups for some of the groups of order 16 also arise as
      direct product factors in some of the automorphism groups of orders
      $16p$ and $16p^2$, which is the reason for their inclusion in this report.

            The groups of order $16p$ and their automorphism groups are
      useful in the determination of groups and automorphism groups of
      groups of higher orders, e.g., $16p^2$, 240, and 192. In \cite{2} we showed
      that the automorphism groups of the groups of orders $8p$ and $8p^2$ displayed
      a very systematic and predictable pattern. A systematic pattern
      analogous to that found for the automorphism groups of groups of
      orders $8p$ and $8p^{2}$ is found in all of the cases we have studied to
      date. The object of this paper is to report the results of the
      calculations of the automorphism groups for the groups of orders $16p$
      and $16p^2$.

            In the presentations given in this paper the modular relation
\begin{equation}
                  x^t \equiv  1\,\, \text{mod}(p)
\end{equation}
      is used to define the quantity $x$. In each case $x$ is a $t$-th root of
      unity. This is not stated explicitly in each case, but this is what
      should be understood when this relation appears. The exponent $t$ is
      usually an even number in the cases given below, so that both $x$ and
      $p-x$ are solutions to this equation. For example, in the presentation
      for the holomorph of $C_{p}$ ($p$ being a prime number $>$ 2), namely
\begin{equation}
            \text{Hol}(C_{p}) : a^p = b^{(p-1)} = a^b*a^x=1,
\end{equation}
      replacing $x$ with $-x$ will yield the same group.

\section {\underline{The groups of order $16p$}}

            The groups of order 16 along with their automorphism groups
      are listed in Table 1. Table 2a gives the groups of order $16p$ with a
      normal sylow $p$-subgroup. Additional material on groups with a normal
      sylow 2-subgroup and those without a normal sylow $p$-subgroup are
      given in Tables 2b and 2c. Table 2a uses the same display format as
      that used in \cite{2} for the groups of order $8p$, namely if $G$ = $C_{p}$ @ $Q$,
      then Aut($G$) = Hol$(C_{p}) \times $ (invariant factor),
      where $Q$ is a group of order 16, $C_{p}$ @ $Q$ is the semidirect product of
      the group $C_{p}$ with $Q$, with $C_{p}$ being the normal subgroup in $G$, and the
      \textquotedblleft invariant factor" is what is given in Table 2a. One can reconstruct
      the groups of order $16p$ using the information given in Tables 1 and
      2a. The following example will show how this is done. Let
\begin{equation}
               G  =  C_{p} @ C_{16}.
\end{equation}
      From Table 1 the presentation for $C_{16}$ is $a^{16}$ =1. The
presentation for $C_{p}$ is just $x^p$ = 1. From Table 2a the group $C_{16}$ acts on
$C_{p}$ as an operator of order 2 by means of the generator $a$ by mapping the generator of order
$p$ into its inverse. Hence the relations for $G$ are:
\begin{equation}
                  x^p=a^{16}=x^a*x=1.
\end{equation}
      All of the other groups in Table 2a of order $16p$ follow this same
      simple pattern. For the case of the groups of order $16p^2$ given
      below when the $p$-group is $C_{p} \times C_{p}$ the operator of order 2 maps one or both
      of the generators of the $p$-group into its inverse. This same format will be used in
      future papers on groups of order 32$p$, $p^3q$, and $p^4q$, among other cases. It
 should be pointed out that the groups of the form $C_{p} \times $ (order 16) are not listed in Tables 2a, 2b, or 2c. The automorphism groups of these groups are just $C_{(p-1)} \times $ Aut(order 16 group).

\section {\underline{The groups of order $16p^2$}}

\subsection      {\underline{  General comments and the $p=3$ case}}

            In \cite{1} and \cite{2} we pointed out certain general relations between the
      automorphism groups of the groups of orders $8p$ and $8p^{2}$. A similar
      situation also occurs in the cases of groups of orders $16p$ and $16p^{2}$.

           The tables giving information on the automorphism groups of the
groups of order $16p^{2}$ are analogous to those found in the report for the order $8p^{2}$ case.

There are 28 direct product groups of the form
\begin{align}
                  (C_{p} \times C_{p}) \times \text{(group of order 16)} \qquad \qquad &\text{(14
cases),} \\
                  C_{p^2}  \times  \text{(group of order 16)} \qquad \qquad &\text{(14
cases).}
\end{align}
      The automorphism groups of these groups are just the direct products
      of the automorphism groups of the $p$-group and the 2-group.
            There are 3 times 28 groups in which the 2-group acts on the $p$-group by means of an operator of order 2:
\begin{align}
      C_{p} @ (\text{2-group}) \times C_{p} \qquad \qquad   & \text{(28 cases),} \qquad & \label{a}\\
      (C_{p} \times C_{p}) @ (\text{2-group})  \qquad \qquad &   \text{(28 cases),} \quad  & \label{b}\\
       C_{p^2} @ (\text{2-group})   \qquad \qquad & \text{(28 cases).}   \quad & \label{c}
\end{align}

      In each case the automorphism group can be determined from the action of the 2-group (read off from Table 2a for the order $16p$ groups). If the $16p^2$ group has the structure (\ref{a}), then its automorphism group is
\begin
{equation}
\text{(factor from Table 2a)} \times \text{Hol}(C_{p}) \times C_{(p-1)}.
\end{equation}
If it has the structure (\ref{b}), then Aut($g$) is
\begin{equation}
\text{(factor in Table 2a)} \times \text{Hol}(C_{p} \times  C_{p})
\end{equation}
and for the form (\ref{c}) we have
\begin{equation}
\text{(factor in Table 2a)} \times \text{Hol}(C_{p^2}) .
\end{equation}
            The groups listed in Table 3a give information on the \textquotedblleft basic groups" of order $16p^2$ and their automorphism groups when the action of the 2-group on the $p$-group is of order 4 or larger. By \textquotedblleft basic groups" we mean groups that have recurrences for all primes $p> 3$. (The one exception in Table 3a is the one coming from $<-2,4|2>$ [\# 13], which occurs in the orders $p\equiv  1$ or 3 mod(8).) Since the automorphism group of $C_{p^2} $ is $C_{p-1} \times C_{p}$, all of the automorphism groups with a $C_{p^2}$ normal $p$-subgroup can be read from the entries in Table 2a. The entries in Table 3a refer to the case when the normal $p$-group is ($C_{p} \times C_{p}$). \\

            Many of the entries in Table 3a make specific reference to the case of $p= 3$. The generalization to primes larger than three is somewhat messy, or cumbersome to list in the table itself, e.g., entry [576](54,4) for one of the ($C_{2} \times C_{2}$) images in ($C_{8} \times C_{2}$). In this case the automorphism group for the case of $p= 3$ has order 576, 54 classes and a center of order 4. The presentations for this and other groups in Table 3a are given below in section 3.3.\\

            The meaning of the entries in Table 3a can best be explained by considering the case of $C_{8} \times C_{2}$. There are two groups of order $16p^{2}$ arising from a ($C_{2} \times C_{2}$) action of ($C_{8} \times C_{2}$) on ($C_{p} \times C_{p}$).\\

      Two more cases arise from a $C_{4}$ action of ($C_{8} \times C_{2}$)  on ($C_{p} \times C_{p}$),
and one additional case for a $C_{8}$ action ($p=3$).\\

The presentations for these groups are
\begin{quotation}
\begin{equation}
a^8=b^2=(a,b)=c^p=d^p=(c,d)= \notag
\end{equation}
 \subsubsection{\underline{The $C_2 \times  C_2$ cases}}
\begin{align}
         & c^a*c=(b,c)=(a,d)=b^d*b=1, &  [a,b]\quad \text{case}, \\
            &c^a*c=(a,d)=c^b*c=d^b*d=1,&     [ab,b]\quad \text{case}.
\end{align}

 \subsubsection{\underline{The $C_4$ cases}}
\begin{align}
& c^a*d^{-1}=d^a*c=(b,c)=(b,d)=1,&  [a]\quad \text{case,}\\
& c^a*d^{-1}=d^a*c=c^b*c=d^b*d=1,&  [ab] \quad\text{case.}
\end{align}
 \subsubsection{\underline{The $C_8$ cases}}
\begin{equation}
 c^a*d^{-1}=d^a*c*d=(b,c)=(b,d)=1   \quad\text{($p=3$ case).}
\end{equation}
\end{quotation}
For the $C_2 \times C_2$ cases, the automorphism groups are:
\begin{align}
(a,b)\text{ type} &\qquad C_{2} \times C_{2} \times C_{2} \times \text{Hol}(C_{p}) \times \text{Hol}(C_{p}), \\
(ab,b)\text{ type}& \qquad [576](54,4),
\end{align}
where [576](54,4) denotes a group of order 576 with 54 conjugacy classes and a center of order 4;
see section 3.3.2 for the presentation.\\

For the $C_{4}$ cases these automorphism groups (for $p=3$) are the same and
isomorphic to
\begin{equation}
D_4 \times  [144],
\end{equation}
where [144] means the complete group of order 144:
\begin{equation}\label{CG144}
      C_3 \times C_3  @ <-2,4|2>.
\end{equation}
The group $<-2,4|2>$ is also known as the quasi-dihedral group of order 16 (or  $QD_{8}$). In the general case ($p>3$) this group goes over into the automorphism group of the group $(C_{p} \times C_{p}) @ C_{4}$. \\

      The groups Aut[($C_{p} \times C_{p}) @ C_{4}$] are complete groups for $p\not\equiv 1$ mod(4). When $p \equiv 1$ mod(4) instead of one group with a $C_{4}$ image, as shown above for the case of $p= 3$, we have four nonisomorphic groups; see Table 4a. As is readily apparent from Table 4a the automorphism groups for these different cases all display the same form, as a function of the prime $p$.\\

             The general expression for the  $(C_{p} \times C_{p}) @ C_{16}$ group's automorphism group when $C_{16}$ acts on the $p$-group as an operator of order 8 is known and is given below (see section 3.2.5 for the case of $p\equiv  1  \bmod (8)$, and section 3.3.1 for $p\not\equiv  1 \bmod (8)$).\\

            General expressions for the $D_4$ and $Q_{2}$ image cases were more difficult to determine. See Table 3b and section 3.3.7 below for more details on this problem.

The group
\begin{equation}
                  [ C_{3} \times C_{3} ] @ <-2,4|2>  \qquad \text{[order 16 action]}
\end{equation}
does have a recurrence in higher orders [$p\equiv  3 \bmod (8)$], but in general this group is not complete, although its automorphism group is a complete group.

      For $p= 3$ there are several additional groups of order $16p^{2}$. These groups are given in Table 5a (groups with a normal sylow 2-subgroup) and Table 5b (groups with no normal sylow subgroup).

      All of the groups in Table 3a have recurrences for primes $p > 3$. For $p > 3$ additional cases arise in much the same manner as do those for $p\equiv  1 \bmod (4)$, or $p\equiv  1 \bmod (8)$,\ldots, etc. in the orders $8p$, $16p$,\ldots, etc.\\

\subsection {\underline{Special cases for groups with $p > 3$}}

            Table 6a, taken from the thesis of R. Nyhlen \cite{5}, gives the
      number of groups of order $16p^{2}$, for various choices of the prime $p$.

The differences between the number of groups for the cases $p = 3$, 5 and 7 and $p\equiv  3$ mod(8), 5 mod(8) and 7 mod(8) are due to the presence of groups with normal sylow 2-subgroups, or groups without a normal sylow subgroup in the cases when $p= 3$, 5 or 7, but which do not arise for larger primes. Table A2 in Appendix 1 gives the subgroups of orders 16 and 32 for some sylow 2-subgroups of $GL(2,p)$. The behavior of the automorphism groups of $(C_{p} \times C_{p}) @ C_{16}$, $(C_{p}
\times C_{p}) @ D_{8}$, and $(C_{p} \times C_{p}) @ Q_{4}$ for the primes $p\equiv  7$ mod(8) as a function
of $p$, shown below, is related to the presence or absence of certain subgroups of orders 16 or 32 in $GL(2,p)$.

There are 166 basic types which appear in all orders. The additional group(s) for $p\equiv  3$ and 7 mod(8) are:
\begin{gather}
\text{$p \equiv  3$ mod(8) (one from} <-2,4|2>\,\, = QD_{8}),\\
\text{$p \equiv  7$ mod(8) (one each from}\,\, C_{16},\,\, D_{8}\,\, \text{and}\,\, Q_{4}).
\end{gather}

Additional groups arising from the groups of order 16 for higher primes are given in Tables 5c and 6b.

For primes $p\equiv  1$ mod(16), the total number of groups according to Nyhlen \cite{5}
(resp. according to R. Laue, see Table 6a) is \\

                  166 + 53 + 24 + 14 (resp. 15) = 257 (resp. 258). \\

      In the discussion of the groups of order $16p$, we noted that if $p$ were equal to 1 mod(4) or 1 mod(8) or 1 mod(16) we obtained additional groups of order $16p$ for these primes. In that discussion the number of groups appearing in order $16p$ for $p\equiv  1$ mod(16) was equal to the number of groups coming from a $C_{4}$ quotient (or the number of new groups for primes $p\equiv  1$ mod(4) only) plus the number of groups coming from a $C_{8}$ quotient (or the number of new groups for primes $p\equiv  1$ mod(8) only) and the number of groups coming from $C_{16}$ (or for those primes $p\equiv  1$ mod(16)). The breakdown Nyhlen uses for these cases is $p\equiv  5$ mod(8) instead of $p\equiv  1$ mod(4) (but not equal to $p\equiv  1$ mod(8)), $p\equiv  9$ mod(16) (instead of $p\equiv  1$ mod(8) only) and $p\equiv  1$ mod(16). We shall follow Nyhlen's usage, but this alternate view makes it clear why the number of groups for the primes $p\equiv  1$ mod(16) is the sum of the number of groups that make their first appearance for $p\equiv 1$ mod(4), and others that make their first appearance in order $p\equiv  1$ mod(8) and finally those other cases that only appear in the case of $p\equiv  1$ mod(16).

Note that the cases arising for $p\equiv  3$ mod(8) and 7 mod(8) do not contribute groups to the case of
$p\equiv 1$ mod(16), like the $p\equiv  5$ mod(8) and 9 mod(16) cases do. We have not independently verified that either number of $C_{16}$ cases listed above is correct; they are taken from Nyhlen's 1919 thesis \cite{5}. (See the section on $p\equiv  1$ mod(16) below for additional information on the groups $(C_{17} \times C_{17}) @ C_{16}$.)\\

\subsubsection {\underline{ Presentations for the actions of some order $16$ groups on $C_{p} \times C_{p}$}}

            The actions of a group $G$ on the group ($C_{p} \times C_{p})$ can be read off from a matrix representation of its generating set when $G$ is represented as a subgroup of Aut$(C_{p} \times C_{p}) = GL(2,p)$. As an example consider the case of a $D_{8}$ action on ($C_{p} \times
C_{p})$. Using the matrix representation given below we have for this group (  $(C_{p} \times C_{p}) @ D_{8}$ ):
\begin{equation}
\begin{split}
            a^p&=b^p=(a,b)=c^8=d^2=c^d*c \\
                  &=a^c*a^{-x}*b^{-x}=b^c*a^x*b^{-x}=(a,d)=b^d*b=1 .
\end{split}
\end{equation}
The other cases follow in a similar manner.

In the automorphism groups that follow, the action of the order 16 group on the group ($C_{p} \times C_{p}$) will be given by the following matrix representations as indicated above. \\

\begin{enumerate}
               \item{$D_8$ case: $p \equiv  \pm 1$ mod(8),\quad $c^8=d^2=c^d*c=1,$}\\
\begin{equation}
c =
    \left(
      \begin{matrix}
           x & x \\
          -x & x
      \end{matrix}
  \right), \qquad
d = \left(
       \begin{matrix}
           1 & 0 \\
           0 & -1
        \end{matrix}
      \right), \qquad
    2x^2 \equiv  1 \,\,\text{mod}(p);
\end{equation}\\

  \item{$<-2,4|2>\,\, = QD_8$ case: $p\equiv  1$ or 3 mod(8),\quad $c^8=d^2=c^d*c^5=1$},

\begin{equation}
c =
    \left(
      \begin{matrix}
           x & x \\
          -x & x
      \end{matrix}
  \right), \qquad
d = \left(
       \begin{matrix}
           1 & 0 \\
           0 & -1
        \end{matrix}
      \right),
  \qquad 2x^2 \equiv  -1 \,\, \text{mod}(p);
\end{equation}\\

\item{$Q_4$ case: $p\equiv  \pm 1$ mod(8),}

\begin{equation}
c =
    \left(
      \begin{matrix}
           x & x \\
          -x & x
      \end{matrix}
  \right), \qquad
d = \left(
       \begin{matrix}
           a & b \\
           b & a
        \end{matrix}
      \right), \qquad
    2x^2 \equiv  1\,\, \text{mod}(p), \qquad
    a^2 + b^2 \equiv  -1\,\, \text{mod}(p);
\end{equation}\\

 \item {$C_{16}$ case: $p\equiv  7$ mod(8), or 9 mod(16),}
                    \end{enumerate}
\begin{equation}
            a^p=b^p=(a,b)=c^{16}=a^c*b^{-1}=b^c*a^{-x}*b^{-y}=1 .
\end{equation}
\underline{$C_{16}: p\equiv  7$ mod(8).}

\begin{equation}\label{AF2:Z}
c = \left(
       \begin{matrix}
           0 & 1 \\
           x & y
        \end{matrix}
      \right).
\end{equation}

Some values for the lower primes are:
\begin{equation}
\text{($p,x,y$)} = (7,1,3),\quad (23,1,4),\quad (31,-1,5),\quad (47,-1,3).
\end{equation}
      Nyhlen \cite{5} gives the following relations for $x$ and $y$:
\begin{equation}\label{AF2:1}
 x = 1 \quad \text{ and} \quad (y^2 + 2)^2 \equiv  2\,\, \text{mod}(p).
\end{equation}
      These relations do not generate matrices of order 16 in the group $GL(2,p)$ for $p= 31$, 47, or 79 but are correct for $p=7$, 23, 71, and 103. The previous relation for $y$ given in (\ref{AF2:1}) seems to work only when the sylow 2-subgroup of $GL(2,p)$ has order 32. The other cases (for $p=31$, 47 and 79) seem to obey the relation
\begin{equation}\label{AF2:2}
x = -1 \quad \text{and} \quad (y^2 - 2)^2 \equiv  2\,\, \text{mod}(p).
\end{equation}
      The group $QD_{16}$ arises below in connection with the automorphism groups of certain groups of order $16p^{2}$. For the case of $C_{16}$ given by (\ref{AF2:1}) a matrix of order 4 that yields the group $QD_{16}$ is

\begin{equation}
d' = \left(
       \begin{matrix}
           0 & 1 \\
          -1 & 0
        \end{matrix}
      \right).
\end{equation}\\

      The presentation for $QD_{16}$ determined by $<c,d'>$ is:\\
\begin{equation}
              (c*d')^2=(c*d'^{-1} )^2=d'^4=c^4*d'*c^{-4}*d'=1.
\end{equation}
      For the most common presentation for $QD_{16}$, namely,
\begin{equation}
                     c^{16}=d^2=c^d*c^{-7}=1,
\end{equation}
      the order 2 matrix $d$ takes the form:
\begin{equation}
d  = \left(
       \begin{matrix}
           1 & 0 \\
           y & -1
        \end{matrix}
      \right)
\end{equation}
with $y$ as given above. For those primes $p$ for which the order 16 matrix is generated by (\ref{AF2:2}), the group $QD_{16}$ does not occur as a subgroup of $GL(2,p)$ (see Table A2). Therefore there are no corresponding matrices $d$ and/or $d'$ here that will generate the group $QD_{16}$ for these primes.\\

\underline{$C_{16}: p \equiv  9$ mod(16).}\\

For $p\equiv  9$ mod(16) we have a matrix of the form ({\ref{AF2:Z}) with
\begin{equation}
\text{($p,x,y$)} = (41,3,0)
\end{equation}
or, in general, $y$ = 0 and $x$ obeys the relation:
\begin{equation}
                        x^8 \equiv  1\,\, \text{mod}(p).
\end{equation}

\subsubsection {\underline{ The automorphism group for the $p\equiv  3$ mod$(8)$ case}}

The one additional\\ group
 occurring in the sequence of groups with $p\equiv  3$ mod(8) has for its automorphism group the group
\begin{equation}
            (C_{p} \times C_{p}) @ ( C_q \times <-2,4|2>).
\end{equation}
      Here $q = (p-1)/2$ and the action of $C_q$ on the $p$-group takes the form
\begin{equation}
C_q  = \left(
       \begin{matrix}
           z & 0 \\
           0 & z
        \end{matrix}
      \right) \qquad
\text{with} \quad z^q \equiv  1\,\, \text{mod}(p).
\end{equation}
  \subsubsection {\underline{ The automorphism groups for the $p\equiv  7$ mod$(8)$ case}}

            The calculations for\\ additional groups arising for the $p\equiv  7$ mod(8) and 9 mod(16) cases are summarized in Table 9. Some discussion of these additional groups is given below.

The three additional groups for the cases with $p\equiv  7$ mod(8) have automorphism groups of the following form:

\underline{$C_{16}$: Action on $(C_{p} \times C_{p})$ is $C_{16}$.}
\begin{gather}
            (C_{p} \times C_{p}) @ (C_3 \times QD_{16}) \quad    \text{for $p= 7$},\\
            (C_{p} \times C_{p}) @ (C_{11} \times (C_3 @ QD_{16}))\quad   \text{for $p= 23$},\\
            (C_{p} \times C_{p}) @ (C_3 \times C_5 \times [128]) \quad \text{ for $p= 31$},
\end{gather}
where [128] is the sylow 2-subgroup of $GL(2,31)\,\, [ = QD_{64}$], and the action is full. The group $C_3 @ QD_{16}$ has the presentation
\begin{equation}
                         c^3=a^{16}=b^2=a^b*a^{-7}=(a,c)=c^b*c=1
\end{equation}
and is isomorphic to $D_{48}$.\\

            For cases larger than $p= 31$ we have not been able to determine the structure of these automorphism groups. The cases for $p= 47$ and higher have automorphism groups that are too large to be determined at this site. The interesting features  that we have not been able to check out in higher orders are the occurrence of the factor $(C_{3} @ QD_{16})$ in the $p= 23$ case and whether the sylow 2-subgroup of $GL(2,p)$ appears in higher orders as it does here for the case of $p= 31$. From the other cases run to date one might have expected to see just the direct product $C_{3} \times  QD_{16}$ here (or more generally the factor [ $C_{q} \times QD_{16}$ ], where $q = (p-1)/2$, for the quotient group).\\

          The orders for this sequence of automorphism groups appear to be $2p^2*(p+1)*(p-1)$.\\

\underline{$D_{8}$ and $Q_{4}$: Actions on $C_p \times C_p.$}

\begin{align}
      (C_{p} \times C_{p}) @ (C_q \times QD_{16}) & \quad \text{for} \,\, D_{8}\,\, \text{and}  \,\, Q_{4} \,\, \text{with $p$\,\,= 7, 23}, \\
(C_{p} \times C_{p}) @ (C_q \times D_{16})&\quad  \text{for} \,\, D_{8}\,\,  \text{and $p= 31$}, \\
(C_{p} \times C_{p}) @ (C_q \times Q_{8})&\quad   \text{for} \,\, Q_{4}\,\, \text{and $p= 31$}.
\end{align}

          The action of the $C_q$  ($q = (p-1)/2)$ on the $C_{p}$'s is given by
\begin{equation}
            a^p=b^p=(a,b)=c^q=a^c*a^{-z}=b^c*b^{-z}=\cdots,\text{where}\,\, z^q \equiv  1\,\,\text{mod}(p).
\end{equation}
            From Table A2 in Appendix 1 we see that $QD_{16}$ is not a subgroup of $GL(2,31)$, and so  cannot appear in the automorphism group of $( C_{31} \times C_{31}) @ D_{8}$ (or $Q_{4}$) in the same way as it does for the cases of $p= 7$ and 23. From Table A2 one might conjecture that if $p$ is equal to $-1$ mod(8), but not equal to $-1$ mod(16), then we get an automorphism group with a structure like that for the $p= 7$ and $p= 23$ cases; otherwise, we get groups with a structure like that for the $p= 31$ case.

  \subsubsection {\underline{The automorphism groups for the $p\equiv  5$ mod$(8)$\! case}}

      The cases in which the action is by means of an operator of order 4 are described above in section 3.1, along with how to obtain the automorphism groups from Table 4a.\\

            Automorphism groups for the $C_{4} \times C_{2}$ image cases are given in Table 4b.\\

  \subsubsection {\underline{The automorphism groups for the $p\equiv  9$ mod$(16)$ case}}

\begin{enumerate}
\item  \underline{The $C_8$ actions}

            For the cases of $C_{8} \times C_{2}$ and $C_{16}$ these groups are essentially groups that result from a $C_{8}$ action on the $p$-groups. The groups in question arising from the $C_{8} \times C_{2}$ group are just the order $8p^{2}$ groups cross a $C_{2}$. The groups in this order coming from $C_{16}$ have the same action as those in the $8p^{2}$ order case with the $C_{8}$ replaced by a $C_{16}$. In each case ($C_{8} \times C_{2}$ and $C_{16}$) the automorphism groups of these groups are just $C_{2} \times  $(automorphism group of the corresponding ($C_{p} \times C_{p}) @ C_{8}$ group) with the same $C_{8}$ action.\\

\item  \underline{The $C_{16}$ action case}

      There is one group of this type and the CAYLEY program has not been able to obtain the automorphism group at this site. The authors believe the automorphism group here is going to be isomorphic to Hol$(C_{p}) wr C_{2}$. \\

      This is the type of behavior also found in other cases in which the $C_{(p-1)}$ action on the $p$-group is off-diagonal:
\begin{equation}
C_{(p-1)}  = \left(
       \begin{matrix}
           0 & a \\
           b & 0
        \end{matrix}
      \right) \qquad
(p,a,b) = (5,1,-1), . . .\\
\end{equation}
\item \underline{Order 16 group actions from $C_8 \times  C_2$}

            The presentations for these groups, taken from R. Nyhlen's thesis \cite{5}, are given in Table 7. Table 7 also gives the automorphism groups for these groups. The case of $p=17$ has been explicitly calculated, those for $p > 17$ are conjectures based upon what happens in the other cases.\\

            The last two groups listed in Table 7 have the same class and order structure and as such cannot be distinguished here. \\

\item \underline{Order 16 group actions from $D_8, QD_8,$ and $Q_4$}
               \end{enumerate}

            The discussion that follows is summarized in Table 9.

            The automorphism groups of
\begin{equation}
                 (C_{p} \times C_{p}) @ D_{8}  \quad \text{and} \quad (C_{p} \times C_{p}) @
Q_{4}
\end{equation}
for $p= 17$ have order 73,984. These automorphism groups are not expressible as direct products and are isomorphic. This automorphism group is a normal subgroup of the group Hol$(C_{17}) wr C_{2}$. A permutation representation of degree 34 for this automorphism group is:
\begin{align}
           a &=  (1,2,3,4,5,6,7,8,9,10,11,12,13,14,15,16,17); \notag \\
           b &=  (2,4,10,11,14,6,16,12,17,15,9,8,5,13,3,7)  \notag \\
             &\qquad (19,21,27,28,31,23,33,29,34,32,26,25,22,30,20,24); \\
           c &=  (2,10,14,16,17,9,5,3)(4,11,6,12,15,8,13,7); \notag \\
           d &=  (1,18)(2,19)(3,20)(4,21)(5,22)(6,23)(7,24)(8,25)(9,26) \notag \\
            & \qquad (10,27)(11,28)(12,29)(13,30)(14,31)(15,32)(16,33)(17,34); \notag
\end{align}
the corresponding presentation is:
\begin{equation}
\begin{split}
a^{17}&=b^{16}=a^b*a^{-3}=b^2*(c*d)^{-2}=a^c*a^8=c^8=d^2\\
&=(a*d)^2*(a^{-1}*d)^2=(b,c)=(b,d)=(c*d)^2*(c^{-1}*d)^2=1.
\end{split}
\end{equation}

            One conjecture as to what the presentation might be for this class of automorphism groups is
\begin{gather}
a^p=b^q=a^b*a^{-x}=b^y*(c*d)^{-2}=a^c*a^t \notag \\
=c^8=d^2=(a*d)^2*(a^{-1}*d)^2\\
=(b,c)=(b,d)=(c*d)^2*(c^{-1}*d)^2=1, \notag
\end{gather}
where $q = (p-1)$, $x$ given by $x^q \equiv  1$ mod($p$) is a primitive root of unity, $y=q/2$, and $t^8 \equiv  1$ mod($p$). \\

            The structure for the automorphism groups of
\begin{equation}
(C_{p} \times C_{p}) @ D_{8} \quad \text{and} \quad (C_{p} \times C_{p}) @ Q_{4}
\end{equation}
for $p= 17$ is:
\begin{equation}
(C_{17} \times C_{17}) @ (C_{16}\text{Y}D_{16}).
\end{equation}
      The amalgamated product $C_{16}$Y$D_{16}$ of the groups $C_{16}$ and $D_{16}$ has the presentation:
\begin{equation}\label{AF2:3}
a^{16}=b^{16}=c^2=(a,b)=(a,c)=b^c*b=a^8*b^{-8}=1,
\end{equation}
and the action on the group $C_{17} \times C_{17}$ is full. For the general case the conjectured structure is:
\begin{equation}
(C_{p} \times C_{p}) @ ( C_{p-1} \text{Y} D_{16} ) \text{ [ for $p \equiv  1$ mod(16) ],}
\end{equation}
and the above presentation (\ref{AF2:3}) is slightly changed, i.e., with $a^{16}$ going into $a^{(p-1)}$ and $a^{8}$ going into $a^{(p-1)/2}$. The action of $C_{p-1}$Y$D_{16}$ on $(C_{p} \times C_{p})$ for $p= 17$ can be read off from the matrix representation of this group:
\begin{equation}
A  = \left(
       \begin{matrix}
           x & 0 \\
           0 & x
        \end{matrix}
      \right), \qquad
B  = \left(
       \begin{matrix}
           6 & 1 \\
           1 & 6
        \end{matrix}
      \right), \qquad
C  = \left(
       \begin{matrix}
           1 & 0 \\
           0 & 16
        \end{matrix}
      \right).
\end{equation}

Here $x$ is given above, and for $p= 17$, $x = 3$. For $p> 17$, an appropriate matrix representation for $D_{16}$ has not been determined.\\

Note also the similarity between this automorphism group and the one given below for the group $(C_{17} \times C_{17}) @ QD_{8}$.\\

The automorphism group of
\begin{equation}
                 (C_{p} \times C_{p}) @ QD_{8}
\end{equation}
(for $p= 17$) has order 36,992, has 65 classes, and is  not expressible as a direct product. This order 36,992 group appears to be
\begin{equation}
            (C_{17} \times C_{17}) @ [ C_{16}\text{Y}(C_{4} wr C_{2})],
\end{equation}
which is the automorphism group for $( C_{17} \times  C_{17}) @ D_4$. For the cases of $p\equiv  3$ mod(8),  Aut[($C_{p} \times C_{p}) @ QD_{8}$] is a complete group, but for $p= 17$ this is not the case; here the automorphism tower goes as follows:
\begin{equation}
(C_{p} \times C_{p}) @ QD_{8}--->[36,992]--->[73,984]--->\text{Hol}(C_{17}) wr C_{2}.
\end{equation}
A permutation representation for this automorphism group is:
\begin{align}
             a&=(1,2,3,4,5,6,7,8,9,10,11,12,13,14,15,16,17); \notag\\
             b&=(2,4,10,11,14,6,16,12,17,15,9,8,5,13,3,7) \notag \\
               &\quad  (19,21,27,28,31,23,33,29,34,32,26,25,22,30,20,24);  \\
             c&=(2,14,17,5)(3,10,16,9)(4,6,15,13)(7,11,12,8); \notag \\
             d&=(1,18)(2,19)(3,20)(4,21)(5,22)(6,23)(7,24)(8,25)(9,26)\notag \\
              &\quad  (10,27)(11,28)(12,29)(13,30)(14,31)(15,32)(16,33)(17,34). \notag
\end{align}
The corresponding presentation is:
\begin{equation}
\begin{split}
                    a^{17}&=b^{16}=a^b*a^{-3}=b^4*(c*d)^{-2}=a^c*a^4\\
                          &=c^4=d^2=(b,c)=(b,d)=(a*d)^2*(a^{-1}*d)^2=1.
\end{split}
\end{equation}

One conjecture as to what the presentation might be for this class of automorphism groups is \\
\begin{equation}
\begin{split}
            a^p&=b^q=a^b*a^{-x}=b^y*(c*d)^{-2}=a^c*a^t\\
               &=c^4=d^2=(b,c)=(b,d)=(a*d)^2*(a^{-1}*d)^2=1,
\end{split}
\end{equation}
where $q = (p-1)$,  $x^q \equiv  1$ mod($p$), $y=q/4$, and $t^4 \equiv  1$ mod($p$). (Is this the same as $(C_{p} \times C_{p}) @ ( C_{p-1}$Y$(C_{4} wr C_{2})$ )  ?)

\subsubsection {\underline{ The $p\equiv  1$ mod$(16)$ case}}

      The relations
\begin{equation}
                a^{17}=b^{17}=(a,b)=c^{16}=a^c*a^{-3}=b^c*b^x=1
\end{equation}
for $x=1,2,3,\ldots,16$ all yield groups of order $16*17^2$. Several cases have duplicate class structure and automorphism groups, so they are not readily distinguishable, using the methods we employed in this study. The number of conjugacy classes and the order of these automorphism groups are given in Table 8.\\

\subsection {\underline{Presentations and comments on various cases}}

\subsubsection {\underline{  The $(C_p \times  C_p) @ C_{16}$   [$C_8$ action]  cases}}
The structure of the group (order of element, number of elements of each order, and the number of conjugacy classes) of order 288 arising as an automorphism group of $(C_3 \times  C_3) @ C_{16}$ is shown in Table 10.
The center of this group is $C_{2}$. The general structure for the automorphism groups of the groups:
\begin{equation}
                        (C_{p} \times  C_{p}) @ C_{16}
\end{equation}
with a $C_{8}$ action is:
\begin{equation}
[(C_{p} \times C_{p}) @ C_t] @ C_{4}  \quad  [t =(p^2 -1),\quad p \not\equiv  1\,\,\text{mod(8)}].
\end{equation}
      The structure ($C_{p} \times  C_{p}) @ C_{t}$ ($t =(p^2 -1) $) arises in connection with the automorphism groups for which the \textquotedblleft invariant factor" is half the expected order. See the discussion of the automorphism groups of the groups of order $p^3*q^2$ in \cite{1} for more details.\\

            The relations for this automorphism group for the $p=3$ case are:
\begin{equation}
\begin{split}
      a^3&=b^3=(a,b)=c^8=a^c*b=b^c*a*b=d^4 \\
         &=a^d*a*b^{-1}=b^d*b^{-1}*a^{-1}=c^d*a*c^{-3}=1;
\end{split}
\end{equation}
and for the general case we have ($p \not\equiv  1$ mod(8)):
\begin{equation}\label{AF2:4}
\begin{split}
      a^p&=b^p=c^t=a^c*b=b^c*a^x*(b^{x^2})=d^4\\
         &=(a,d)=b^d*(c*a^x*c^{-1})^{-1} = c^d*c^{(-p)}=1;
\end{split}
\end{equation}
here $t=p^2-1$, and $x$ is the exponent in the relation for the holomorph of $C_{p}$:
\begin{gather}
            a^p=b^{(p-1)}=a^b*a^x=1; \quad \text{i.e.,}\quad x^{(p-1)} \equiv  1\,\, \text{mod}(p);\\
\text{where} \quad (p,x)=(3,1), (5,3), (7,2), (11,5), (13,2), (17,6),\ldots .\notag
\end{gather}

Note in the $p= 3$ case we have $c^d*a*c^{-3}$, whereas in the expression for the general case the $a$ does not appear. Both relations give rise to the same group for $p= 3$.\\

For $p= 17$ the group (\ref{AF2:4}) is not the automorphism group of any group of the form $(C_{17} \times  C_{17}) @ C_{16}$ with either a $C_{8}$ or a $C_{16}$ action. The order of this group is 332,928, with a trivial center. It might be of interest to see if this group with ($p,x$) = (17,6), which has 170 conjugacy classes, is a complete group. The automorphism group of this group ($p = 17$ case) is too large to run at this site. The automorphism groups of all groups of the form $(C_{17} \times  C_{17}) @ C_{16}$ with a
$C_{8}$ action have a nontrivial center, $C_{34}$ for the case $C_{17} @ C_{16} \times
C_{17}$ or $C_{2}$ for the other cases.  \\

\subsubsection{\underline{ Aut = $576$ with ncl = $54$ and Z = $C_2 \times  C_2$ groups.
        $(C_2 \times  C_2)$ action cases}}

       The group in the $C_{2} \times C_{2}$ image with order 576, 54 classes
and center $C_{2} \times C_{2}$ (which is \# 5394 in the Small Group Library (SGL) \cite{13})
has the representation:
\begin{equation}
\begin{split}
a&=(1,2,3); \quad b=(2,3); \quad  c=(4,5); \\
 d&=(1,6)(2,7)(3,8)(4,9)(5,10)(11,12,13,14),
\end{split}
\end{equation}
which looks like a modification of a wreath product. If the 4-cycle in $d$ is omitted, this is just the wreath product $(S_3 \times C_{2}) wr C_{2}$. For general primes $p$ just replace the Hol($C_3$) by Hol($C_{p}$); e.g., for $p=5$:
\begin{equation}
\begin{split}
      &a---->(1,2,3,15,16), \quad   b---->(2,3,16,15),\\
      &d---->(1,6)(2,7)(3,8)(4,9)(5,10)(11,12,13,14)(15,17)(16,18).
\end{split}
\end{equation}
      A presentation for this series of groups is:
\begin{equation}
\begin{split}
                      a^p&=b^q=a^b*a^x=c^2=(a,c)=(b,c)=d^4=(a,d^2)\\
                     &=(b,d^2)=(c,d^2)=(a*d)^2*(a^{-1}*d^{-1})^2\\
                     &=a*d*b*d^{-1}*a^{-1}*d^{-1}*b^{-1}*d\\
                     &=a*d*c*d^{-1}*a^{-1}*d^{-1}*c*d\\
                     &=(b*d)^2*((d*b)^{-1})^2\\
                     &=b*d*c*d^{-1}*b^{-1}*d^{-1}*c*d\\
                     &=(c*d)^2*(c*d^{-1})^2=1,
\end{split}
\end{equation}
where $q=p-1$, $x^q\equiv 1\,\,\text{mod}(p)$, and $x$ takes on the values
\begin{equation}
                  (p,x) = (3,1), (5,2), (7,2),\ldots .
\end{equation}

   \subsubsection {\underline{Aut = $2304$ cases $(C_2 \times  C_2)$ actions}}

       The group of order 2304 with 90 classes and a center $C_{2}$ has the following general form:
\begin{equation}
[ \text{Hol}(C_{p}) \times \text{Hol}(C_{p}) \times (C_{2} \times C_{2}) wr C_{2} ] @ C_{2},
\end{equation}
where the $C_{2}$ acts on the holomorphs as in the wreath product action, and $C_{2}$ acting on the group of order 32 produces Hol($C_{4} \times C_{2}$). The presentation for this set is:
\begin{equation}
\begin{split}
                  a^2&=b^2=c^2=(a,b)=(a*c)^4=(a*c*b*c)^2\\
                     &=(b*c)^4=d^2=(a,d)=b^d*b*a=c^d*a*c*a\\
                     &=e^p=f^q=h^p=k^q=e^f*e^x=h^k*h^x\\
                     &=(e,h)=(e,k)=(f,h)=(f,k)=e^d*h\\
                     &=f^d*h*k^{-1}=h^d*e=k^d*e*f^{-1}\\
                      &=(a,e)=(b,e)=(c,e)=(a,f)=(b,f)=(c,f)\\
                     &=(a,h)=(b,h)=(c,h)=(a,k)=(b,k)=(c,k)=1.
\end{split}
\end{equation}
Here $q=p-1$, and $x$ is the exponent in the relation for the holomorph: e.g., $(p,x)=(3,1),$ (5,2), (7,3), \ldots . The group generated by ($a$,$b$,$c$,$d$) is Hol($C_{4} \times C_{2}$). The group generated by
($d$,$e$,$f$,$h$,$k$) is Hol($C_{p}) wr C_{2}$. \\

\subsubsection {\underline{The order $1152$ automorphism groups\! $(3$ cases$)$\! $(C_2 \times  C_2)$ \!actions}}
\!There are three different groups of order 1152 arising as automorphism groups of the groups of order 144. These groups have different numbers of conjugacy classes and take the form (in the general $16p^{2}$ cases):
\begin{align}
      [1152](90,4) \textup{ SGL } \# 90180 :  \quad &[\text{Hol}(C_{p}) \times C_{2} \times C_{2}] wr C_{2}, \\
      [1152](63,2) \textup{ SGL } \# 97975 : \quad &[\text{Hol}(C_{p}) \times \text{Hol}(C_{p}) \times D_4 \times C_{2}] @ C_{2},\\
      [1152](81,2) \textup{ SGL } \# 138373: \quad &[\text{Hol}(C_{p}) \times \text{Hol}(C_{p}) \times D_4 \times C_{2}] @ C_{2}.
\end{align}

In the [1152](63,2) case the action of the $C_{2}$ on $D_4 \times C_{2}$ produces the group Hol($C_{8}$), while in the third case the group $D_4$Y$D_4$ arises. In the last two cases, the action of the $C_{2}$ on the holomorphs is as in the case of a wreath product. Presentations for these cases are: \\

\underline{[1152](63,2):}
\begin{equation}
\begin{split}
                            a^4&=b^2=a^b*a=c^2=(a,c)=(b,c)=d^2\\
                               &=(a,d)=b^d*c*a^{-1}*b=c^d*c*a^2\\
                               &=e^p=f^q=h^p=k^q=e^f*e^x=h^k*h^x=d^2\\
                               &=(e,h)=(e,k)=(f,h)=(f,k)\\
                               &=e^d*h=f^d*h*k^{-1}=h^d*e=k^d*e*f^{-1}\\
                               &=(a,e)=(b,e)=(c,e)=(a,f)=(b,f)=(c,f)\\
                               &=(a,h)=(b,h)=(c,h)=(a,k)=(b,k)=(c,k)=1.
\end{split}
\end{equation}
Here $q$ and $x$ are as in the 2304 case above. The group generated by ($a$,$b$,$c$,$d$) is Hol($C_{8}$).\\

\underline{[1152](81,2) case:}
\begin{equation}
\begin{split}
                           a^4&=b^2=a^b*a=c^2=(a,c)=(b,c)=d^2\\
                              &=a^d*a=(b,d)=c^d*c*a^2\\
                              &=e^p=f^q=h^p=k^q=e^f*e^x=h^k*h^x=d^2\\
                              &=(e,h)=(e,k)=(f,h)=(f,k)\\
                              &=e^d*h=f^d*h*k^{-1}=h^d*e=k^d*e*f^{-1}\\
                              &=(a,e)=(b,e)=(c,e)=(a,f)=(b,f)=(c,f)\\
                              &=(a,h)=(b,h)=(c,h)=(a,k)=(b,k)=(c,k)=1.
\end{split}
\end{equation}
Here the $q$ and $x$ are as in the 2304 case above. The group generated by ($a$,$b$,$c$,$d$) in this case is $D_4$Y$D_4$. In the last two cases the group generated by ($d$,$e$,$f$,$h$,$k$) is the wreath product:
Hol($C_{p}) wr C_{2}$.

\subsubsection {\underline{ Aut($S_3 \times  S_3 \times  C_2 \times  C_2$), order $6912$ case}}

The automorphism group of $S_3 \times  S_3 \times C_{2} \times C_{2}$ of order [6912] and its recurrences has a similar form, namely:
\begin{equation}
[\text{Hol}(C_{p}) \times \text{Hol}(C_{p}) \times ((C_{2} \times C_{2}) wr C_{2}) @ C_3 ] @ C_{2},
\end{equation}
where the action of the $C_{2}$ on the holomorphs is as in the wreath product case, and the action of $C_{2}$ on the order 96 group produces the group of order 192 without a normal sylow subgroup with 14
classes and a trivial center. See \# 54 in Table 4.6b of \cite{1} or Table 4e of \cite{12}. The presentation for this set of automorphism groups is:
\begin{equation}
\begin{split}
                  a^2&=b^2=c^2=(a,b)=((a,c),a)=((b,c),a)\\
                     &=((b,c),b)=x^3=a^x*b=b^x*b*a=(c,x)\\
                     &=d^2=(a,d)=b^d*b*a=(c,d)=x^d*x\\
                    &=e^p=f^q=h^p=k^q=e^f*e^x=h^k*h^x=d^2\\
                     &=(e,h)=(e,k)=(f,h)=(f,k)\\
                     &=e^d*h=f^d*h*k^{-1}=h^d*e=k^d*e*f^{-1}\\
                     &=(a,e)=(b,e)=(c,e)=(a,f)=(b,f)=(c,f)\\
                     &=(a,h)=(b,h)=(c,h)=(a,k)=(b,k)=(c,k)\\
                     &=(x,e)=(x,f)=(x,h)=(x,k)=1.
\end{split}
\end{equation}
Note. For these ($C_{2} \times C_{2}$) extensions, the automorphism groups given here can also be put into the form:
\begin{equation}
                        [\text{Hol}(C_{p}) \times \text{Hol}(C_{p}) ] @ G[32]
\end{equation}
or $G$[64], $G$[128] or $G$[192] respectively with the same presentations. In effect a certain group is acting on the holomorphs as in the wreath product. We believe the representation of the automorphism groups given above is the easiest one(s) to visualize what is happening. \\

  \subsubsection {\underline{ $(C_p \times  C_p) @ Q_2 \times  C_2$ case}}

The order 1728 group arising for the other $Q_{2}$ image has an automorphism group tower that terminates after two iterations at order 6912. For $p=5$, the corresponding group has order 9600. Its automorphism tower also terminates with two iterations at the order 38,400, [$9600\rightarrow 19200\rightarrow 38400$]. For the higher-order recurrences, the orders of the automorphism groups have a curious pattern:
\begin{equation}
 |\text{Aut}(g)| = 2^6*p^2 * (3,6,9,15,18)\quad \text{for}\ p=(3,5,7,11,13)\ \text{respectively}.
\end{equation}
The automorphism towers for the cases $p > 5$ have not been followed to completion. We believe the towers terminate in the third iteration. For the case of $p=7$ we have the sequence:
\begin{equation}
 784\longrightarrow  28,224\longrightarrow  56,448\longrightarrow  112,896.
\end{equation}
We have not been able to determine whether the last group in this sequence is complete, or whether this chain continues to higher orders.\\

The automorphism groups here are related to the automorphism groups of [$(C_{p} \times C_{p}) @ Q_{2}$]. The automorphism groups of these groups have the structure:
\begin{equation}
            (C_{p} \times C_{p}) \text{@(group of order 24}(p-1)).
\end{equation}
See Appendix 2 for details on these automorphism groups. The automorphism group for the group ($C_{p} \times C_{p}) @ Q_{2} \times C_{2}$  has the structure:
\begin{equation}
(C_{p} \times C_{p} \times C_{2} \times C_{2}) \text{@(group of order 24}(p-1)).
\end{equation}
The action of the group of order 24($p-1$) on the $p$-subgroup is the same in both cases. The group of order 24($p-1$) acts on the group ($C_{2} \times C_{2}$) as follows:\\

$p \equiv  3$ mod(8) cases:
\begin{equation}
\begin{split}
                 f^2&=h^2=(f,h)\\
                    &=a^3=b^2=((a*b)^2*a^{-1}*b)^2=e^q=(a,e)=(b,e)=  (***)\\
                    &=f^a*h=h^a*f*h=(b,f)=h^b*f*h=(e,f)=(e,h)=1,\\
                    &\quad\quad\quad  \text{where $q = (p-1)/2$.}
\end{split}
\end{equation}
For $p\equiv  5$ mod(8) just replace the relations in the line $(***)$ by those for the $p\equiv  5$ mod(8) case of Aut[($C_{p} \times C_{p}) @ Q_{2}$], and similarly for the other cases (see Appendix 2a and
Table A3). Note that in the relations for Aut[$(C_{p} \times C_{p}) @ Q_{2}$] given in Appendix 2 for the $p\equiv  1$ mod(8) and 7 mod(8) cases the order 3 and the order $2^n$ elements are interchanged; so in the third line above, the generators $a$ and $b$ should be interchanged.

\subsection{\underline{The automorphism groups of the groups of order
$16p^{2}$ with a $D_4$ or}\\
\underline{a $Q_{2}$ action on ($C_{p} \times C_{p}$)}}

For a given $p$ the automorphism groups for the groups with a $D_4$ action for $p= 3,7,11,19$ and 23 have the same class and order structure. The case of $p\equiv 5$ mod 8 ($p = 5, 13$,\ldots ) is somewhat different in that the automorphism groups coming from the groups $D_4 \times C_{2}$ and $C_{4} @ C_{4}$ (both $D_4$ and the $Q_{2}$ image case) are isomorphic and distinct from the automorphism groups coming from the groups $D_{8}$ and $Q_{4}$ with a $D_4$ action on the $p$-group. The case of $p=17$ appears to
follow that for $p\equiv  5$ mod(8) in that the $D_4 \times C_{2}$, $C_{4} @ C_{4}$ (both $D_4$ and $Q_{2}$
images) have the same automorphism groups. The reason for this behavior follows from the structure for these automorphism groups, namely that
\begin{equation}
\text{Aut}[ (C_{p} \times C_{p}) @ H ] = [ (C_{p} \times C_{p}) @ T \times (C_{2} \times C_{2}) ] @ C_{2}.
\end{equation}
      If the group $H$ acts on ($C_{p} \times C_{p}$) as $D_4$, then $T$ is isomorphic to the group $D_4$Y $C_{p-1}$. If $H$ acts as $Q_{2}$, then $T$ is isomorphic to $Q_{2}$Y$C_{p-1}$ . For the cases of $p\equiv  1$ mod(4), $D_4$Y$C_{p-1} \cong  Q_{2}$Y$C_{p-1}$. The automorphism groups here appear to split into several distinct types with the following general presentations. This pattern is also found in the automorphism groups of the groups of order $32p^2$ when the action of the order 32 group on the $p$-group ($C_{p} \times C_{p}$) is by a $D_4$ action \cite{6}.\\

\subsubsection {\underline{$p\equiv  3$ mod($8$)}}
\begin{itemize}
\item  \underline{$D_{4}$ image case:}
\begin{equation}
\begin{split}
                  a^8&=b^2=a^b*a^5=c^p=d^p=(c,d)\\
                     &=c^a*c^{(-x)}*d^x=d^a*c^{(-x)}*d^{(-x)}\\
                     &=(b,c)=d^b*d\\
                     &=e^2=f^y=(a,f)=(b,f)=(e,f)=c^f*c^{(-t)}=d^f*d^{(-t)}\\
                     &=a^2*e*a^{-2}*e=(b,e)=(c,e)=(d,e)\\
                     &=(a*e)^2*(a^{-1}*e)^2=1,
\end{split}
\end{equation}
            \text{where ($p$,$x$,$y$,$t$) = (3,1,1,1), (11,4,5,3), (19,3,9,4),....}
\item      \underline{$Q_{2}$ image case}:\\
\begin{equation}
\begin{split}
                  a^8&=b^2=a^b*a^5=c^p=d^p=(c,d)\\
                     &=c^a*c^{(-x)}*d^x=d^a*c^{(-x)}*d^{(-x)}\\
                     &=(b,c)=d^b*d\\
                     &=e^2=(b*e)^4=a^2*e*a^{-2}*e\\
                     &=(a*e)^2*(a^{-1}*e)^2=(c,e)=(d,e)\\
                     &=a*e*a^{-1}*b*e*b\\
                     &=f^y=(a,f)=(b,f)=(e,f)=c^f*c^{(-t)}=d^f*d^{(-t)}=1,
\end{split}
\end{equation}
            \text{where ($p$,$x$,$y$,$t$) = (3,1,1,1), (11,4,5,3),....}
\end{itemize}
\subsubsection {\underline{$p\equiv  5$ mod($8$)}}
\begin{itemize}
\item \underline{$D_4$ [from $D_4 \times C_{2}$ and $C_{4} @ C_{4}$] and $Q_{2}$
image cases.}
\end{itemize}
The following presentation will yield the automorphism groups for this set of groups:
\begin{equation}
\begin{split}
                        a^p&=a^b*a^x=d^q=a^d*a^y\\
                           &=b^d*a*b^{-1}*a \qquad \qquad      (*)\\
                           &=(c,d)=e^2=(a,e)=(c,e)=(d,e)\\
                           &=c^2=b^4=b^2*e*b^{-2}*e\\
                           &=a*c*a*c*a^{-1}*c*a^{-1}*c\\
                           &=a*c*b*c*a^{-1}*c*b^{-1}*c\\
                           &=b*c*b*c*b^{-1}*c*b^{-1}*c\\
                           &=b*e*b*e*b^{-1}*e*b^{-1}*e\\
                           &=(b*e*b^{-1}*c)^2=1,
\end{split}
\end{equation}
where ($p,x,q,y$) = (13,$-5$,3,$-3$), (29,12,7,13),... and
$q=(p-1)/4$, $x^4 \equiv  1$ mod($p$), and $y^q \equiv  1$ mod($p$).\\

For the case of $p=5$, one should delete relations involving the generator $d$ to get a correct presentation. The particularly offending word is the one indicated above by ($*$). If this is deleted, then setting $p=5$,... and $d=1$, one can generate the $p=5$ member of this sequence. \\

Note if we set $e=1$, then we get the corresponding sequence of automorphism groups for the groups of order $8p^{2}$ with the $D_4$ acting on the corresponding $p$-group. \\

      A degree 30 permutation representation for the case with $p=13$ is:
\begin{equation}
\begin{split}
            a&=(1,2,3,4,5,6,7,8,9,10,11,12,13);\\
            b&=(1,8,12,5)(2,3,11,10)(4,6,9,7)(27,29)(28,30);\\
            c&=(1,14)(2,15)(3,16)(4,17)(5,18)(6,19)(7,20)(8,21)\\
              &\qquad \qquad(9,22)(10,23)(11,24)(12,25)(13,26);\\
            d&=(2,10,4)(3,6,7)(5,11,13)(8,12,9)\\
                  &\qquad \qquad (15,23,17)(16,19,20)(18,24,26)(21,25,22);\\
            e&=(27,28).
\end{split}
\end{equation}
\begin{itemize}
\item \underline{$D_4$ actions: Cases from $D_{8}$ and $Q_{4}$.}
\end{itemize}
The structure of these automorphism groups is
\begin{equation}
                  [(C_{p} \times C_{p}) @ D_4 \times C_{2} \times C_{2} ] @ C_{4}
\end{equation}
or
\begin{equation}
                  [ D_p wr C_{2} \times C_{2} \times C_{2}] @ C_{4}.
\end{equation}
The actions of the $C_{4}$ on $D_p wr C_{2}$ turn this group into $(C_{p} @ C_{4}) wr C_{2}$, and the actions on $C_{2} \times C_{2}$ turn this group into the group $<2,2,2>$. An alternate representation for this automorphism group's structure is:
\begin{equation}
[(C_{p} wr C_{2}) \times C_{2}] @ (C_{8} \times C_x), \quad \text{where}\,\, x = (p-1)/4.
\end{equation}
The presentations and relations given below are most easily viewed using this decomposition for this series of groups.\\

A permutation representation for the case of $p=5$ with degree 18 is:
\begin{equation}
\begin{split}
                  a&=(1,2,3,4,5);\\
                  b&=(2,3,5,4)(11,12,13,14,15,16,17,18);\\
                  c&=(1,6)(2,7)(3,8)(4,9)(5,10);\\
                  d&=(12,16)(14,18).
\end{split}
\end{equation}
A prescription for the general case can easily be obtained for a permutation representation for any finite $p\equiv  5$ mod(8). For the case of $p=13$ we have
\begin{equation}
\begin{split}
           a&=(1,2,3,4,5,6,7,8,9,10,11,12,13);\\
            b&=(2,6,13,9)(3,11,12,4)(5,8,10,7)\\
                  &\quad (27,28,29,30,31,32,33,34);\\
            c&=(1,14)(2,15)(3,16)(4,17)(5,18)(6,19)(7,20)\\
                  &\quad (8,21)(9,22)(10,23)(11,24)(12,25)(13,26);\\
            d&=(28,32)(30,34);\\
            e&=(2,4,10)(3,7,6)(5,13,11)(8,9,12)\\
                  &\quad(15,17,23)(16,20,19)(18,26,24)(21,22,25).
\end{split}
\end{equation}
In presentation form we have
\begin{equation}
\begin{split}
            a^p&=b^8=a^b*a^x=c^2=d^2=(a,d)=b^d*b^3=(c,d)\\
               &=(a*c)^2*(a^{-1}*c)^2=a*c*b*c*a^{-1}*c*b^{-1}*c\\
               &=(b*c)^2*(b^{-1}*c)^2=(b*d*b^{-1}*c)^2\\
               &=e^3=(b,e)=(d,e)=a^e*a^(-y)=(c,e)=1;
\end{split}
\end{equation}
here $x$ and $y$ are determined by the following modular relations:
\begin{equation}
x^4 \equiv  1\,\,\text{mod}(p) \quad \qquad \text{and} \qquad y^3 \equiv  1\,\,\text{mod}(p).
\end{equation}

\subsubsection {\underline{$p\equiv  7$ mod($8$)}} The $D_4$ and $Q_2$ image cases.\\
\begin{itemize}
\item \underline{$D_4$ image:}
\end{itemize}
\begin{equation}
\begin{split}
                  a^8&=b^2=a^b*a=c^p=d^p=(c,d)\\
                      &=ca*c^{-x}*d^x=d^a*c^{-x}*d^{-x}\\
                     &=(b,c)=d^b*d\\
                     &=e^2=f^y=(a,f)=(b,f)=(e,f)=c^f*c^{-t}=d^f*d^{-t}\\
                     &=a^2*e*a^{-2}*e=(b,e)=(c,e)=(d,e)\\
                     &=(a*e)^2*(a^{-1}*e)^2=1,
\end{split}
\end{equation}
where ($p,x,y,t$) = (7,2,3,4),....
\pagebreak\
\begin{itemize}
\item \underline{$Q_{2}$ image:}
\end{itemize}
\begin{equation}
\begin{split}
                  a^8&=b^4=a^4*b^{-2}=a^b*a=c^p=d^p=(c,d)\\
                     &=c^a*c^{-x}*d^x=d^a*c^{-x}*d^{-x}\\
                     &=cb*c^{x1}*d^{x2}=d^b*c^{x2}*d^{-x1}\\
                     &=e^2=f^y=(a,f)=(b,f)=(e,f)=c^f*c^{-t}=d^f*d^{-t}\\
                     &=a^2*e*a^{-2}*e=(b,e)=(c,e)=(d,e)\\
                     &=(a*e)^2*(a^{-1}*e)^2=1,\\
\end{split}
\end{equation}
where ($p,x,x1,x2,y,t$) = (7,2,4,2,3,4),....

\subsubsection {\underline{$p\equiv  1$ mod($16$): here only one case is available, $p= 17$}}

      In this order the automorphism groups coming from the $D_4$ and $Q_{2}$ actions of $C_{4} @ C_{4}$ and $D_4 \times C_{2}$ on the group $C_{17} \times C_{17}$ appear to yield the same automorphism group. For the $D_4$ image from the $C_{4} @ C_{4}$ case, CAYLEY gives the following presentation for the automorphism group:
\begin{equation}
\begin{split}
            a^2&=d^2=(a,b)=(a,e)=(c,e)=(d,e)=b^3*e^{-1}*b^{-1}*e\\
               &=b^2*e^2*b^{-1}*e^{-2}=(a*c)^2*(a*c^{-1})^2\\
               &=(a*d*b^{-1}*d)^2=a*c*a*e^2*c^{-3}=a*c^2*a*d*c^{-2}*d\\
               &=a*c^2*d*a*d*c^{-2}=b*c*b*c^{-1}*b^{-1}*c*b^{-1}*c^{-1}\\
               &=b*c*d*c^{-1}*b^{-1}*c*d*c^{-1}=(b*d)^2*(b^{-1}*d)^2\\
               &=c^4*e^4=(c*d)^2*(c^{-1}*d)^2=b^2*c^2*d*b*d*c^{-2}\\
               &=b^2*c^{-1}*b^{-1}*c*d*b^2*d=1.
\end{split}
\end{equation}

      The general case for $D_4 \times C_{2}$ and $C_{4} @ C_{4}$ may be given by the following presentations:
\begin{equation}
\begin{split}
            a^p&=b^q=a^b*a^x=c^4=d^2=(a,c)=(b,c)=(b,d)\\
               &=(a*d)^2*(a^{-1}*d)^2=b^t*d^{-1}*c^{-1}*d^{-1}*c^{-1}\\
               &=e^2=(a,e)=(b,e)=(d,e)\\
               &=e*c^2*e*c^{-2}=(e*c)^2*(e*c^{-1})^2=1,
\end{split}
\end{equation}
where $q=p-1$, $x^q\equiv 1$\,\,mod($p$), and $t=q/4$.

This set of possible relations has not been checked out for $p > 17$. The orders for the groups are too large for CAYLEY to run here. The group $<a,b,c,d>$ is the automorphism group for the corresponding $8*p^2$ cases (verified for $p= 41$). \\

The other case is for the automorphism groups with a $D_4$ image from the groups $D_{8}$ and $Q_{4}$. The case of $D_{8}$ yields, from CAYLEY, the relations for Aut($g$):
\begin{equation}
\begin{split}
            b^2&=e^2=(a,b)=(a,d)=(a,e)=(b,c)=a*b*a*d*b*d\\
               &=a*c^3*a^{-1}*c^{-1}=a^2*c*a^{-2}*c^{-2}\\
               &=c*d*(c*e)^2*d^{-1}=b*d^2*b*e*d^{-2}*e\\
               &=b*d^2*e*b*e*d^{-2}=(b*e)^4=c^2*d^{-1}*c*e*c^{-1}*e*d\\
               &=(d*e)^2*(d^{-1}*e)^2=a^2*e*d^{-1}*e*d^{-5}=1.
\end{split}
\end{equation}

For $p > 17$ the orders of the groups and their automorphism groups are too large to be handled here. The next cases in this grouping occur for $p= 97$, 113 and 129.\\

\subsubsection {\underline{$p\equiv  9$ mod($16$)}}

The first nontrivial case occurs at $p=41$. For the cases of $p\equiv  1$ mod($2^n$) with $n >$ 1 and  $p >5$ (including the case of $p\equiv  1$ mod(8)) we conjecture that: \\
\begin{quotation}
A presentation of the automorphism groups for the $D_4 \times C_{2}$ (with a $D_4$ action) and $C_{4} @ C_{4}$ (with either a $D_4$ or a $Q_{2}$ action) extension cases is given by the following:
\begin{gather*}
            a^p=b^q=a^b*a^x=c^4=d^2=(a,c)=(b,c)=(b,d)\\
            =(a*d)^2*(a^{-1}*d)^2=b^t*d^{-1}*c^{-1}*d^{-1}*c^{-1}\\
            =e^2=(a,e)=(b,e)=(d,e)\\
            =e*c^2*e*c^{-2}=(e*c)^2*(e*c^{-1})^2=1,\\
            \text{with $q$, $x$, and $t$ as in the $p\equiv  1$ mod(16) case,}\\
 \text{but here with $p\equiv  1$ mod($2^n$).}
\end{gather*}
\end{quotation}

We have not been able to compare these conjectured relations (for $p\equiv  1$ mod(8) and $p > 17$) to ones calculated by CAYLEY, because the group orders are too large for CAYLEY to obtain the corresponding
automorphism groups in question at this site. \\

3.4.6. \underline{Some general comments on the automorphism groups arising from a}\\
\underline{group of order 16 with a $D_4$ or a $Q_{2}$ action on
($C_{p} \times C_{p}$).}\\

In each of the above cases the group generated by $<a,b,c,d>$
is the group $(C_{p} \times C_{p})$ @ [2-group], where the action of the 2-group on the $p$-group listed in Table 3b ($<-2,4|2>$, $D_{8}$ or $Q_{4}$) is full; i.e., no normal subgroup of the 2-group commutes with the $p$-group. \\

A more instructive way to see how these automorphism groups are generated, and differ from one another, is to look at permutation representations for some of these groups: \\

\underline{$D_4$ image ($p= 3$ case):}
\begin{equation}
\begin{split}
                  a&=(2,9,4,6,3,7,5,8)(10,12)(11,13),\\
                  b&=(4,5)(6,9)(7,8),\\
                  c&=(10,11),\\
                  d&=(1,3,2)(4,7,6)(5,8,9),\\
                  e&=(1,5,4)(2,9,6)(3,8,7).
\end{split}
\end{equation}
\underline{      $Q_{2}$ image  ($p = 3$ case):}
\begin{equation}
\begin{split}
                  a&=(2,9,4,6,3,7,5,8)(10,12)(11,13),\\
                  b&=(4,5)(6,9)(7,8)(10,12)(11,13),\\
                  c&=(10,11),\\
                  d&=(1,3,2)(4,7,6)(5,8,9),\\
                  e&=(1,5,4)(2,9,6)(3,8,7).
\end{split}
\end{equation}

The first 9 \textquotedblleft letters" generate the complete group
\begin{equation}
                  (C_3 \times C_3) @ <-2,4|2>
\end{equation}
of order 144. The \textquotedblleft 10 and 11" part of the group have $<-2,4|2>$ acting as a wreath product on the generator $c$ in the $D_4$ case, and a rather odd similar behavior in the $Q_{2}$ case. For cases with $p >$ 3 the
permutations for the group $(C_{p} \times C_{p}) @ <-2,4|2>$ become more complicated and one needs to add the part in which the cyclic group $C_q$ acts on the group ($C_{p} \times C_{p}$), where $C_q$ is as
in Appendix 2a. The $D_4$ image case above corresponds to \# 8299 of order 576 in the small group
library, and the $Q_2$ image case corresponds to \# 8300 or order 576. \\

For the cases $D_4 \times C_{2}$, or $C_{4} @ C_{4}$, and $p=5$, a permutation representation for the automorphism group is: \\
\begin{equation}
\begin{split}
                     a&=(1,2,3,4,5),\\
                  b&=(6,7),\\
                     c&=(1,3,4,2)(6,13)(7,14),\\
                  d&=(1,8)(2,9)(3,10)(4,11)(5,12).
\end{split}
\end{equation}
This is a variant of a wreath product group. If the factors (6,13)(7,14) were in the generator $d$ instead of in the generator $c$, the group would be
\begin{equation}
                  [\text{Hol}(C_5) \times C_{2}] wr C_{2}.
\end{equation}
The action of $d$ on $C_5 \times C_5$ is the same as in the case of the wreath product, and likewise so is the $C_{2}$ action of $c$ on ($C_{2} \times C_{2}$).\\

For the cases of a $D_4$ action from the groups $D_{8}$ or $Q_{4}$, we have the following permutation representation of the automorphism groups:\\
\begin{equation}
\begin{split}
                  a&=(1,2,3,4,5),\\
                  b&=(12,16)(14,18),\\
                  c&=(2,3,5,4)(11,12,13,14,15,16,17,18),\\
                  d&=(1,6)(2,7)(3,8)(4,9)(5,10).
\end{split}
\end{equation}
It might be of interest to note that the automorphism tower for the $D_4$ action on $C_5 \times C_5$ coming from the $D_{8}$ (or $Q_{4}$) group has, for its next two terms, groups of orders 25,600 and 819,200.\\

\section {\underline{Conclusions and Problems}}

      The program outlined in the introduction has been shown to be feasible for the groups of orders $16p$ and $16p^{2}$; namely, for each group for which the order $16p^{2}$ group can be expressed as a presentation involving the running variable $p$, one can determine a presentation for the corresponding group's automorphism group.\\

The authors have not obtained expressions for all cases; alternate noncomputer methods seem the most appropriate way to proceed to finish off the remaining cases for the orders $16p^{2}$, involving the groups that occur only in orders 1 or 7 mod(8). The remaining problems occur for the automorphism groups of the groups
\begin{equation}
            (C_{p} \times C_{p}) @ D_{8} \quad \text{and} \quad (C_{p} \times C_{p}) @ Q_{4},
\end{equation}
where the action on the $p$-group is by the group $D_4$ and, for primes $p$ of the form $p\equiv  1$ mod(8), the action on the $p$-group is by the full $D_8$ or $Q_4$.\\

The factors given in Tables 2a, 3a, and 4a (i.e., those factors independent of $p$, are an invariant of the 2-group and are determined by the 2-group associated with the extension of the group $G$[16], namely the kernel of the homomorphism:
\begin{equation}
            f: \text{[order 16 group]}---> \text{Aut}(C_{p} \times C_{p}).
\end{equation}
      The groups of order $32p$ offer a better illustration of this effect and will be described more fully in the next paper devoted to the automorphism groups of the groups of order $32p$.\\

      In the above analysis in which the action of the 2-groups on the $p$-group was $C_{2} \times C_{2}$, the authors have the feeling that one could just look at the ($C_{2} \times C_{2}$) actions on the $p$-group and tell if their automorphism groups would be isomorphic. To be more precise, the conjecture is that all of the groups of order $16p^{2}$ (arising from a $C_{2} \times C_{2}$ action) with the same automorphism group could be made to have the same $C_{2} \times C_{2}$ actions on the $p$-group. One should note here that in these cases the 2-groups are not isomorphic. For the cases in which the automorphism groups take the form \\
\begin{equation}
                  \text{Hol}(C_{p}) \times \text{Hol}(C_{p}) \times (\text{some group)},
\end{equation}
the actions are
\begin{equation}
                  a^p=b^p=a^c*a^x=(b,c)=(a,d)=b^d*b^x= \cdots .
\end{equation}
This effect seems reasonable, but in the \textquotedblleft cross-action cases", where at least one of the generators from the 2-group acts on both of the generators of the $p$-group, this may still be the case. In these cross-action cases the automorphism groups are no longer a direct product, but in these cases one may still be able to have the actions of the 2-group on the $p$-group be the same for those groups that yield
isomorphic automorphism groups. The number of cases dealt with here is comparatively few in number to what is supposed to appear in the orders $32p^2$. Therefore, the order $32p^2$ groups might be a better place to test this conjecture if one cannot come up with a more traditional group-theoretic proof of the correctness of this assertion or show that this conjecture is false.\\

      In the case of $p\equiv  7$ mod(8), the automorphism groups for the groups of order 16 ($C_{16}$, $D_{8}$ and $Q_{4}$) acting on ($C_{p} \times C_{p}$) seem to have a structure that depends upon the prime in question; i.e., they seem to behave very differently than for the other primes in which the sylow 2-subgroups of $GL(2,p)$ are isomorphic, for different $p$'s. For the cases of $p\equiv  7$ mod(8), can one find a general structural pattern to fit these cases, e.g., some sort of periodic recurrence pattern for these structures as the prime increases? In this connection we conjecture that for $p\equiv  -1$ mod(8) (but not $p\equiv  -1$ mod(16)) the automorphism groups for
\begin{equation}
      (C_{p} \times C_{p}) @ D_{8} \qquad   \text{and}  \qquad  (C_{p} \times C_{p}) @ Q_{4}
\end{equation}
are isomorphic to each other and to $(C_{p} \times C_{p}) @ ( C_q \times  QD_{8})$.\\

      If $p\equiv  -1$ mod($2^n$), where $n$ is larger than 3, then we have
\begin{equation}
\text{Aut}(C_{p} \times C_{p}) @ D_{8}    =   (C_{p} \times C_{p}) @ ( C_q \times  D_{16})
\end{equation}
and
\begin{equation}
\text{Aut}(C_{p} \times C_{p}) @ Q_{4}    =   (C_{p} \times C_{p}) @ ( C_q \times  Q_{8} ).
\end{equation}
Here $q = (p-1)/2$. We do not have enough cases on hand to make any conjectures as to what happens for the automorphism groups of ($C_{p} \times C_{p}) @ C_{16}$. These are situations in which the orders of the groups and their automorphism groups are becoming very large, and where other more group-theoretic (noncomputer) methods for determinng the automorphism groups of these groups might be more appropriate to use.\\

      In some of the automorphism groups given above the sylow 2-subgroups of $GL(2,p)$ occur, and it would be useful to have explicit matrix representations that correspond to the presentations given in Table A2 of Appendix 1. The usefulness of such matrix representations is that then the actions of these 2-groups on the $p$-group ($C_{p} \times C_{p}$) can be determined very easily. The explicit forms for these matrix representations are given in Tables 11a and  11b. These naturally split into two forms: \\

      \underline{a.) for  $p\equiv  - 1$ mod($2^n$):}
\begin{equation}
a =
    \left(
      \begin{matrix}
           0 & 1 \\
           1 & x
      \end{matrix}
  \right), \qquad
b = \left(
       \begin{matrix}
           1 & 0 \\
           x & -1
        \end{matrix}
      \right)
\end{equation}
corresponding to the presentations
\begin{equation}
               a^{t1} = b^2 = a^b*a^{-t2} = 1,
\end{equation}
where
\begin{equation}
       t1=2^{n+1} \quad  \text{ and} \quad t2=(2^n -1).
\end{equation}
The case of $p\equiv  -1$ mod(8) is given in section 3.2.1. An alternate presentation is given by Carter and Fong in \cite{8}. See Table 11a for more details. \\

      \underline{b.) for $p\equiv  + 1$ mod($2^n$)} we have
\begin{equation}
a =
    \left(
      \begin{matrix}
           1 & 0 \\
           0 & x
      \end{matrix}
  \right), \qquad
b = \left(
       \begin{matrix}
           0 & 1 \\
           1 & 0
        \end{matrix}
      \right)
\end{equation}
with  $x$  being a $2^n$-th root of unity in $GF(p)$, corresponding to the presentation:
\begin{equation}
            a^t = b^2 = (a*b)^2*(a^{-1}*b)^2=1 \quad  (t=2^n).
\end{equation}
      See Table 11b for comments on other wreath products.\\

\section{Appendices}
\subsection{ Sylow 2-subgroups of $GL(2,p)$}

      In attempting to get explicit presentations for the automorphism
      groups in Table 3b, the following observations about the sylow 2-subgroups
      of $GL(2,p)$ were useful:\\

       \underline{$<-2,4|2>$ = $QD_8$}
\begin{quotation} is the sylow 2-subgroup of Aut$(C_p \times C_p)$
                  [$GL(2,p)$ for $p\equiv 3$ mod(8)].
\end{quotation}

       \underline{$C_4 wr C_2$}
\begin{quotation} is the sylow 2-subgroup of Aut($C_p \times C_p$)
                  [$GL(2,p)$ for $p\equiv 5$ mod(8)].
\end{quotation}

       \underline{$D_8$ is NOT the sylow 2-subgroup of Aut($C_7 \times C_7$).}
\begin{quotation}
                  The sylow 2-subgroup of $GL(2,7)$ is \underline{$QD_{16}$}, i.e.,
                  group number 50 of order 32 in the tables of
                  M. Hall and J. K. Senior \cite{7}.
\end{quotation}
The factor in the $p= 17$ case does not have the same order as the sylow 2-subgroup of Aut($C_{17} \times C_{17}$). The factor here, however, is consistent with being the quotient of this sylow 2-subgroup by ($C_2 \times C_2$).\\

The other interesting case(s) of $p= 41$ (9 mod(16)) is too large for CAYLEY to obtain this factor. It might be of interest here to note, however, that the quotient group of the sylow 2-subgroup of $GL(2,41)$ by ($C_2 \times C_2$) is the group $C_4 wr C_2$.

\begin{gather}
            \text{CONJECTURE:}\notag \\
                  \text{Aut}[(C_{41} \times C_{41}) @ (C_4 @ C_4)] =
            [C_{41} \times C_{41} \times C_2 \times C_2] @ [(C_4 wr C_2) \times C_5]
\end{gather}
      For the cases $p\equiv  3$ mod(8), 5 mod(8) and 9 mod(16) the sylow 2-subgroups of $GL(2,p)$ have the orders 16, 32, 128, respectively, but for $p\equiv  7$ mod(8) and 1 mod(16) this is not the case. In these last two cases we have:

\begin{gather}
 \text{\underline{$p \equiv  7$ mod(8)}} \quad (23,2^5),\quad (31,2^7),\quad (47,2^6),\quad (71,2^5),\quad (79,2^6), \\
                  (103,2^5),\quad (127,2^9),\quad (151,2^5),\quad (167,2^5),\quad (191,2^8),\notag\\
             \text{\underline{ $p\equiv  1$ mod(16)}} \quad(17,2^9),\quad (97,2^{11}),\quad (113, 2^9),
\end{gather}
      where ($x,y$) = ($x$ = prime, $y$ = order of sylow 2-subgroup).\\

       The interesting point here is that the number of groups of order $16p^2$ does not change with $p$ in these sequences (i.e., for the cases with primes $p\equiv  7$ mod(8) or $p\equiv  1$ mod(16)) even though the order of the sylow 2-subgroups does change with $p$. The sylow 2-subgroups of these two sequences have a rather interesting subgroup structure in order for this to occur (see Table A2). For the case of $p\equiv  7$ mod(8) and $G$ = $(C_p \times C_p) @ C_{16}$, the automorphism groups of these groups seem to \textquotedblleft fill up the 2-group actions on the $p$-group"; i.e., the order of the automorphism groups seems to always have a factor of $2^n$, where the number $2^n$ is the order of the sylow 2-subgroup of $GL(2,p)$ [ = Aut($C_p \times C_p$)]. We have only been able to verify this for the cases up to $p= 31$. It would be useful to see if other more traditional methods can be brought to bear on this problem. The orders of these groups become very large, the next cases being $p= 47$, 71, and 89. The most interesting primes here would seem to be of the form $p\equiv  -1$ mod($2^n$) for $p > 31$ and
$n > 3$, e.g., $p= 47$, 223, 191 and 127.\footnote{Note that 223 $\equiv -1$ mod($2^5$) but 223 $\not\equiv -1$ mod($2^6$), and 191 $\equiv -1$ mod($2^6$) but 191 $\not\equiv -1$ mod ($2^7$), etc.} A similar thing seems to happen in the $p\equiv  1$ mod(16) cases for which the wreath product Hol($C_p)wr C_2$ or the holomorph of ($C_p \times C_p$) is the automorphism group for the group ($C_p \times C_p) @ C_{16}$.\\

      In Tables 6a and 6b where the number of groups of order $16p^2$ is given as a function of the prime $p$, the values of the primes are expressed as $p\equiv  3$ mod(8), 5 mod(8), etc. From the point of
view of the sylow 2-subgroups of the groups $GL(2,p)$, a better breakdown would seem to be according to the sequence $p\equiv  +1$ mod(4) or $p\equiv  -1$ mod(4) (omitting the case $p= 2$). The sylow 2-subgroup of $GL(2,p)$ for $p\equiv  1$ mod(4), i.e., $p= 2^n + 1$, is isomorphic to $C_y wr C_2$ (where $y=2^n$). If $p\equiv  -1$ mod(4), then $p= 2^n - 1$, and the order of the sylow 2-subgroup is $4 * 2^n$, is isomorphic to the quasi-dihedral group $QD_t$ (where $t= 2^{n+1}$), and has the presentation:
\begin{equation}
      a^{t1}=b^2=a^b*a^{-t2} = 1, \quad \text{where}\,\, t1=2^{n+1},\quad  t2= 2^{n-1}.
\end{equation}
If $p= 7$, then we have $n$ = 3, and the sylow 2-subgroup of $GL(2,p)$ has the presentation:
\begin{equation}
            a^{16} = b^2 = a^b*a^{-7} = 1;
\end{equation}
if $p= 31$, then $n$ = 5 and the sylow 2-subgroup for $GL(2,p)$ has the presentation:
\begin{equation}
            a^{64} = b^2 = a^b*a^{-31} = 1.
\end{equation}

      Therefore, one has the following classification of the sylow 2-subgroups for\linebreak $GL(2,p)$:\\
\underline{if $p\equiv  -1$ mod($2^n$)},
\begin{quotation}
 then the sylow 2-subgroups for all
primes $p$ with the same $n$ are isomorphic, and isomorphic to $QD_t$ (where $t$ = $2^{n+1}$).
 \end{quotation}
\underline{if $p\equiv  +1$ mod($2^n$)},
\begin{quotation}
 then the sylow 2-subgroups for all primes $p$ with the same $n$ are isomorphic, and isomorphic to $C_y wr C_2$ (where $y = 2^n$), or to $C_{2^n} wr C_2$.
\end{quotation}

            The subgroups of order 16 and 32 contained in the sylow-2 subgroups for $GL(2,p)$ for selected primes are shown in Table A2.

\begin{tabular}{|c|c|c|c|c|c|c|} \hline
  \multicolumn{7}{|c|}{Table A1} \\ \hline
  \multicolumn{7}{|c|}{order structure for some sylow 2-subgroups of $GL(2,p)$}
  \\ \hline
prime & \multicolumn{6}{|c|}{class/order structure of $GL(2,p)\dagger $ } \\
\cline{2-7}
  $p$ & order 2 & order 4 & order 8 & order 16& order 32& order 64 \\ \hline
   3  &  (1,4)  &  (2,4)  &    2\^{}2    &      &      &     \\
\hline
   5  & (1,2,4) & (1\^{}2,   &    4\^{}2    &      &      &
\\
           &    &  2\^{}5,8) &           &          &          &
\\ \hline
   7   &  (1,8)  &  (2,8)  &    2\^{}2    &    2\^{}4   &      &
\\ \hline
   17  & (1,2,16)& (1\^{}2,   & (1\^{}4,     & (1\^{}8, 2\^{}92, &   16\^{}8
&      \\
           &         & 2\^{}5,16) & 2\^{}22,16\^{}2)&   16\^{}4)   &       &
       \\ \hline
      23   &  (1,8)  &  (2,8)  &    2\^{}2    &    2\^{}4   &      &
  \\ \hline
      31   &  (1,32) &  (2,32) &    2\^{}2    &    2\^{}4   &   2\^{}8    &
  2\^{}16 \\ \hline
      41   & (1,2,8) & (1\^{}2,   & (1\^{}4,     &    8\^{}4   &      &
    \\
           &         &  2\^{}5,8) & 2\^{}22,8\^{}2) &          &          &
       \\ \hline
      47   &  (1,16) &  (2,16) &    2\^{}2    &    2\^{}4   &   2\^{}8    &
    \\ \hline
      71   &  (1,8)  &  (2,8)  &    2\^{}2    &    2\^{}4   &      &
  \\ \hline
     79   &  (1,16) &  (2,16) &    2\^{}2    &    2\^{}4   &   2\^{}8    &
   \\ \hline
      97   & (1,2,32)& (1\^{}2,   & (1\^{}4,     & (1\^{}8,    & (1\^{}16,
&  32\^{}16    \\
           &         & 2\^{}5,32)  & 2\^{}22,32\^{}2)&
2\^{}92,16\^{}4)&2\^{}376,32\^{}8)&           \\ \hline
     127 $\ast$ & (1,128) & (2,128) &    2\^{}2    &    2\^{}4   &   2\^{}8
   &   2\^{}16 \\ \hline
     191 $\ast\ast$ &  (1,64) &  (2,64) &    2\^{}2    &    2\^{}4   &   2\^{}8
   &   2\^{}16 \\ \hline
\multicolumn{7}{|c|}{} \\
\multicolumn{7}{|c|}{$\ast$ The $p= 127$ case has 32 classes of elements of
order 128, each with    } \\
\multicolumn{7}{|c|}{two elements in each class, and 64 classes of elements
    } \\
\multicolumn{7}{|c|}{of order 256 with two elements in each class.
    } \\
\multicolumn{7}{|c|}{} \\
\multicolumn{7}{|c|}{$\ast\ast$ The $p= 191$ case has 32 classes of elements of
order 128, each with    } \\
\multicolumn{7}{|c|}{two elements in each class. } \\
\multicolumn{7}{|c|}{} \\
\hline
\multicolumn{7}{|c|}{} \\
\multicolumn{7}{|c|}{$\dagger$ Here ($a$\^{}$x$, $b$\^{}$y$, $c$) means that there are $x$ classes of $a$} \\
\multicolumn{7}{|c|}{elements each, $y$ classes of $b$ elements each, and one} \\
\multicolumn{7}{|c|}{class of $c$ elements each. For example, (1\^{}1, 8\^{}1) $\equiv (1,8)$} \\
\multicolumn{7}{|c|}{means there are two classes of elements of order 2, one} \\
\multicolumn{7}{|c|}{with 1 element in it and the other with 8 elements. The form} \\
\multicolumn{7}{|c|}{1\^{}2 means there are two classes each with one element in it, etc.} \\
\multicolumn{7}{|c|}{} \\
 \hline
\end{tabular}

\begin{tabular}{|c|c|c|c|} \hline
  \multicolumn{4}{|c|}{Table A2} \\ \hline
\multicolumn{4}{|c|} {Subgroups of some sylow 2-subgroups of $GL(2,p)$ }\\
\hline
prime & order & subgroups of order 16 & subgroups of order 32 \\
  & of group & & \\ \hline
   3    &    16   &      $ QD_8$               &           none          \\
\hline
   5    &    32   &$ C_4 \times  C_4, C_4$Y$Q_2, <2.2|2> $&   $C_4 wr C_2$
  \\ \hline
   7    &    32   & $      C_{16}, D_8, Q_4     $&   $        QD_{16}$
\\ \hline
   17   &   512   &$ C_{16}, C_8 \times  C_2, D_8, QD_8, Q_4 $&$  C_{32},
C_{16} \times  C_2, C_8 \times  C_4,$ \\
         &         &$ C_4 \times  C_4, C_4$Y$Q_2, <2,2|2> $&$  C_8$Y$Q_2,
C_4 wr C_2, D_{16},$  \\
          &          &                         &$  QD_{16}, Q_8, <4,4|8>$
and   \\
          &          &                         &  number 22 of order 32  \\
\hline
   23   &    32    & $      C_{16}, D_8, Q_4$      &$         QD_{16}  $
\\ \hline
   31   &   128    &$       C_{16}, D_8, Q_4 $&      $  D_{16},  Q_8,
C_{32}  $\\ \hline
   41    &  128  &$ C_{16}, C_8 \times  C_2, D_8, QD_8, Q_4 $&$ C_4 \times
C_8, C_4 wr C_2, \text{and} \#$'s \\
   &        &$ C_4 \times  C_4, C_4$Y$Q_2, <2,2|2> $&   21, 22, 26 of order
32 \\ \hline
     47   &    64    &$       C_{16}, D_8, Q_4         $&$       D_{16},
Q_8,  C_{32} $  \\ \hline
    71   &    32    & $      C_{16}, D_8, Q_4     $ &     $      QD_{16}   $
\\ \hline
    79   &    64   &$       C_{16}, D_8, Q_4   $&$      D_{16},  Q_8,
C_{32} $ \\ \hline
   127  &   512    &$       C_{16}, D_8, Q_4 $   &$      D_{16},  Q_8,
C_{32}$  \\ \hline
  191   &   256    &$       C_{16}, D_8, Q_4        $&$      D_{16},  Q_8,
C_{32}    $ \\ \hline

\multicolumn{4}{|c|}{} \\
\multicolumn{4}{|c|}{ The sylow 2-subgroups of $GL(2,p)$ for $p= 103$, 151, and
167        } \\
\multicolumn{4}{|c|}{are isomorphic to the $p= 7$ case. } \\
\multicolumn{4}{|c|}{} \\
\hline
\multicolumn{4}{|c|}{} \\
\multicolumn{4}{|c|}{The sylow 2-subgroups of $GL(2,p)$ for $p\equiv  1$ mod($2^n$)
are isomorphic to    } \\
\multicolumn{4}{|c|}{} \\
\multicolumn{4}{|c|}{$C_t wr C_2\,\, (t = 2^n)$} \\
  \multicolumn{4}{|c|}{} \\
  \multicolumn{4}{|c|}{with presentations:       } \\
\multicolumn{4}{|c|}{} \\
\multicolumn{4}{|c|}{$ a^t=b^2=(a*b)^2*(a^{-1}*b)^2=1.  $                 }
\\
\multicolumn{4}{|c|}{} \\
\multicolumn{4}{|c|}{The sylow 2-subgroups of $GL(2,p)$ for $p\equiv  3$ or 7 mod(8)
are isomorphic   } \\
  \multicolumn{4}{|c|}{to the quasi-dihedral group of the appropriate order
with the      } \\
\multicolumn{4}{|c|}{presentation:                  } \\
\multicolumn{4}{|c|}{} \\
\multicolumn{4}{|c|}{$ a^x=b^2=a^b*a^{-y}=1, $   where $x = |G|/2$ and $y =
(x/2) - 1.$            } \\
\multicolumn{4}{|c|}{} \\
\multicolumn{4}{|c|}{Here $|G|$ = the order of the sylow 2-subgroup of
$GL(2,p)$.} \\
\multicolumn{4}{|c|}{} \\
\hline
\end{tabular}

\subsection{\underline{ Appendix 2a. Automorphism groups for ($C_p \times C_p) @ Q_2$}(\cite{9})}

            The presentations of these groups of order $8p^2$ can be read off from the representation of $Q_2$ in $GL(2,p)$.\\

            A matrix representation of $Q_2$ for the presentation
\begin{equation}
                             a^4=b^4=a^2*b^2=a^b*a=1
\end{equation}
is
\begin{equation}
a =
    \left(
      \begin{matrix}
           x & y \\
           y & -x
      \end{matrix}
  \right),\qquad
b = \left(
       \begin{matrix}
           0 & 1 \\
          -1 & 0
        \end{matrix}
      \right),
\end{equation}
where $(x^2 + y^2) \equiv  -1 \bmod (p)$.\\

The automorphism groups of these groups have the structure:
\begin{equation}
                  (C_p \times C_p) @ \text{(group of order 24($p-1$)).}
\end{equation}
      The group of order 24($p-1$), in this sequence of automorphism groups, depends upon the prime $p$ as follows:\\

\begin{tabular}{|c|c|c|}
\hline
prime & quotient group & $q$-factor \\ \hline
$p\equiv 3$ mod(8) & $GL(2,3) \times   C_q $& $q=\frac{p-1}{2} $\\ \hline
$p\equiv 5$ mod(8) & $SL(2,3) @ C_4 \times C_q$ &$q=\frac{p-1}{4}$ \\ \hline
$p\equiv 7$ mod(8) & $<2,3,4> \times C_q$ & $q=\frac{p-1}{2}$\\ \hline
$p\equiv 1$ mod(8) & $GL(2,3)@ C_4 \times C_q$ & $(p-1)=q*2^n$ \\ \hline
\end{tabular}
\linebreak\\

      In the cases where  $p\equiv  1$, 3, or 5  mod(8) these groups of order 24($p-1$) can be written in the form $SL(2,3) @ C_{p-1}$ with the presentation:
\begin{equation}
\begin{split}
                  a^2&=b^3=(a*b*a*b*a*b^{-1})^2=(a,c)=(b,c)\\
                     &=c^x=c^y*(a*b)^4=1, \quad\text{where $x=(p-1)$ and $y=x/2$.}
\end{split}
\end{equation}
A matrix representation for this series of groups is given by:
\begin{equation}
a = \left (
     \begin{matrix}
       0 & s \\
       t & 0
     \end{matrix}
    \right), \quad
b = \left(
     \begin{matrix}
             -1&1 \\
              -1 & 0
      \end{matrix}
      \right), \quad
c=\left(
    \begin{matrix}
       z&0 \\
       0&z
      \end{matrix}
    \right),
\end{equation}
where the matrix $c$ is just the center of the group $GL(2,p)$ and $H$ = $<a,b>$ is a representation of $GL(2,3)$ by $2 \times  2$ matrices over the field $GF(p)$. For the $p\equiv  1$ mod(8) case we have for the first few cases:
\begin{align*}
                  p &= (17,41,73,89,97,113,\ldots);\\
                  s &= (2,3,10,\ldots);\\
                  t &= (9,14,22,\ldots);\\
                  z &= (3,6,5,3,5,3,\ldots).
\end{align*}
In the general case we have $s$, $t$ and $z$ being given by the following relations:
\begin{equation}
                  s^8 \equiv  1\,\,\text{mod}(p),\quad t = (p+1)/s, \quad \text{and}
\quad             z^{p-1} \equiv  1\,\,\text{mod}(p)
\end{equation}
(see \cite{10}). In the other case the group $GL(2,3)$ is replaced by the group $<2,3,4>$. In this case we have the following presentation:
\begin{equation}
                  a^{-2}*b^3=a^{-2}*c^4=a^{-1}*b*c=
                    d^x=d^y*a^2=(a,d)=(b,d)=(c,d)=1.\\
\end{equation}
      Explicit presentations for the $Q_2$ image cases are given in Table A3 for the cases of $p\equiv  1$, 3 or 5 mod(8). In these presentations the roles of $a$ and $b$ are interchanged in the $p\equiv  3$ and 5 mod(8) cases.\\

       2.a.1. The case of $p\equiv  7$ mod(8). The automorphism groups in this sequence are given by the following  presentation:
\begin{equation}
\begin{split}
                  a^3&=b^4=(a,b^2)=a*b*(a*(b^{-1}))^3\\
                     &=c^p=d^p=(c,d)=c^a*c*d^{-1}=d^a*c \label{A2:1a} \\
                     &=c^b*c^v*d^x=d^b*c^y*d^w\\
                     &=e^q=(a,e)=(b,e)=c^e*c^{-z}=d^e*d^{-z}=1,
\end{split}
\end{equation}
where
\begin{gather}
(p,v,w,x,y)=(7,1,-1,1,5),\quad (23,1,-1,7,3),\quad (31,1,-1,12,5),\\
\quad(47,1,-1,18,26),\quad (71,1,-1,7,20),\ldots \notag
\end{gather}
 and $q=(p-1)/2$, $z^q \equiv  1$ mod($p$). Other values for these quantities for various primes can be found in Tables A4 and A5 below.\\

            Here the matrices are
\begin{equation}
a = \left (
     \begin{matrix}
       -1 & 1 \\
       -1 & 0
     \end{matrix}
    \right), \quad
b = \left(
     \begin{matrix}
             v&x \\
              y & w
      \end{matrix}
      \right), \quad \label{A2:A2a.2}
\end{equation}
where $a$ and $b$ generate the group $<2,3,4>$.\\

In this form $b$ has order $4*q=4*(p-1)/2$. For many (but not all !) primes of the form $p\equiv  7$ mod(8) one can represent the generator $b$ for the presentation (\ref{A2:1a}) with $v=-w=1$, and $x*y \equiv  -2$ mod($p$). \\

            An alternate representation for this automorphism group is based upon the group $\textless 2,3,4 \textgreater \times  C_q$ and is generated by matrices of the form (\ref{A2:A2a.2}), and is (for $p= 7$)
\begin{equation}
\begin{split}
               a^3&=b^{4q}=(a,b^2)=(a*b^{-1})^4*b^{-(p-5)}\\
                  &=c^p=d^p=(c,d)=c^a*c*d^{-1}=d^a*c \label{A2:A2a.3a}  \\
                  &=c^b*c^{-1}*d^{-x}=d^b*d*c^{-y}\\
                  &=(a*b)^4*b^2=1,
\end{split}
\end{equation}
and for $p > 7$ we have
\begin{equation}
\begin{split}
               a^3&=b^{4q}=(a,b^2)=(a*b^{-1})^4*b^{-(p-5)} \\
                  &=c^p=d^p=(c,d)=c^a*c*d^{-1}=d^a*c  \qquad      \label{A2:A2a.3b}  \\
                  &=c^b*c^{-1}*d^{-x}=d^b*d*c^{-y} \\
                  &=(a*b)^2*a^{-1}*b^{-1}*(a*b^{-1})^2*a^{-1}*b=1.
\end{split}
\end{equation}

One can also find representations for $b$ for the forms (\ref{A2:A2a.3a}) and (\ref{A2:A2a.3b}) with $v=-w=1$. We do not have closed form expressions for the values of $x$ and $y$ for either form of the above presentations. A collection of values for various primes $p$ is given below in Table A6. Note these values of $x$ and $y$ are different than those used in the presentation found in (\ref{A2:1a}).\\

            This is a problem in the representation theory of $\textless2,3,4\textgreater$ in terms of 2 $\times $ 2 matrices over the field $GF(p)$. See Appendix 2b below for details.\\

\begin{tabular}{|c|c|c|}
\hline
\multicolumn{3}{|c|}{Table A3}\\ \hline
\multicolumn{3}{|c|}{automorphism groups for $(C_p \times  C_p) @ Q_2$}\\
\hline
primes & presentations & modular relations \\ \hline
   3 mod(8)  &$ a^3=b^2=((a*b)^2*a^{-1}*b)^2= $ & $y = (p-1)/2$    \\
             &$ c^p=d^p=(c,d)=c^a*c*d^{-1}=d^a*c=  $&$ x^{p-1} \equiv 1$ mod($p$) \\
             &$ c^b*c^{-1}*d^{-x}=d^b*d=   $           &$ y^q \equiv  1$ mod($p$) \\
                   & $e^q=c^e*c^{-y}=d^e*d^{-y}= $& \\
                    & $ (a,e)=(b,e)=1                   $&        \\ \hline
         5 mod(8)  &$ a^3=(a,b^2)=(a*b)^4=(a*b^{-1})^4=  $&   $q=(p-1)/4$
\\
                    &$ c^p=d^p=(c,d)=c^a*c*d^{-1}=d^a*c=     $  & $ x^4 \equiv  1$
mod($p$)  \\
                   &$ c^b*d^{-x}=d^b*c^{-1}=                   $ &$  y^q \equiv  1 $
mod($p$) \\
                   &$     e^q=(a,e)=(b,e)=c^e*c^{-y}=     $  &
   \\
                    & $        d^e*d^{-y}=1             $  &
  \\ \hline
         7 mod(8)   &    see discussion below           &             \\
\hline
1 mod(8) &$ a^2=b^3=(a*b*a*b*a*b^{-1})^2= $         &   $q=(p-1)$      \\
                  & $   (a,c)=(b,c)=c^q=c^t*(a*b)^4=     $  &   $t=q/2$     \\
                   & $   d^p=e^p=(d,e)=d^a*e^{-x}=         $  &
\\
                  &$   e^a*d^{-y}=d^b*d*e^{-1}=e^b*d=     $   &$ x^8  \equiv  1$
mod($p$)\\
                  & $      d^c*d^{-z}=e^c*e^{-z}=1       $  & $y = (p+1)/s$
  \\
                    &                                 &$ z^{(p-1)} \equiv 1$
mod($p$)\\ \hline
\end{tabular} \\

\subsection{Appendix 2b. Question on matrix representations of
            the Coxeter group $\textless 2,3,4\textgreater$ by 2 $ \times $ 2 matrices
over $GF(p)$}

      The Coxeter group comes up in connection with getting the presentations for the automorphism groups of the groups
\begin{equation}
                         (C_p \times C_p) @ Q_2 \quad \text{when $p\equiv  7$ mod(8).}
\end{equation}
      For many cases the following matrices will give a matrix representation for the Coxeter group $\textless2,3,4\textgreater$ :

\begin{equation}
a =
    \left(
      \begin{matrix}
           -1 & 1 \\
           -1 &  0
      \end{matrix}
  \right), \qquad
b = \left(
       \begin{matrix}
           1 & x \\
           y & -1
        \end{matrix}
      \right).
\qquad \qquad \label{A2:A2b.1}
\end{equation}\\

      Here the matrix $a$ has order 3 and $b$ is of order 4. The presentation for this matrix representation is
\begin{equation}
            a^3=b^4=(a,b^2)=a*b*(a*(b^{-1}))^3=1 \qquad \qquad \label{A2:A2b.2}\\
\end{equation}
when the entries in the matrix $b$ obey the relation:
\begin{equation}
                  x*y \equiv  -2\,\,\text{mod}(p).
\end{equation}
      Problem: Find an algebraic method for determining the values of $x$ and $y$ in the above matrix $b$. The following are the values found by trial and error (actually a computer run for various primes).\\

\begin{tabular}{|c|c|c|c|c|c|} \hline
\multicolumn{6}{|c|}{Table A4} \\ \hline
\multicolumn{6}{|c|}{ Matrix elements for $<2,3,4>$ in (\ref{A2:A2b.1})}\\ \hline
prime & ($x,y$) & prime & ($x,y$) & prime &($x,y$) \\
\hline
             7   &   (1,5)  &   103   &   none   &   223  &   (57,43)  \\
                 &   (2,6)  &         &          &        & (101,117)  \\
            23   &   (7,3)  &   127   &   none   &        & (106,122)  \\
                 &  (20,16) &         &          &        & (180,166)  \\
            31   &   (12,5) &   151   &   none   &   239  &  (31,131)  \\
                 &  (26,11) &         &          &        &  (46,187)  \\
            47   &  (18,26) &   167   & (54,68)  &        &  (52,193)  \\
                 &  (21,29) &         & (99,113) &        & (108,208)  \\
                 &  (20,14) &   191   &  (8,143) &   263  &    none    \\
                 &  (33,27) &         & (42,100) &        &            \\
            71   &   (7,20) &         & (48,183) &   271  &  (114,19)  \\
                 &  (51,64) &         & (91,149) &        &  (209,35)  \\
            79   &  (19,29) &   199   &  (12,33) &        &  (236,62)  \\
                 &  (50,60) &         & (166,187)&        &  (252,157) \\
\hline
\end{tabular}\\

      Matrix representations of the form (\ref{A2:A2b.1}) do not exist for the primes 103, 127, 151 and 263. The values listed in the above table are apparently the only possible values for $x$ and $y$ yielding the above representation.\\

    Alternate forms are required for the primes 103, 127, 151, and 263 and probably others larger than 263. A few possible choices for the matrix $b$, which together with $a$ above obeys the relations
(\ref{A2:A2b.2}), are:

\begin{equation}
b = \left(
       \begin{matrix}
           v & x \\
           y & w
        \end{matrix}
      \right),
\end{equation}\\
      where ($v,x,y,w$) are given in the following table:\\

\begin{tabular}{|c|c|}  \hline
\multicolumn{2}{|c|}{Table A5}\\ \hline
\multicolumn{2}{|c|}{ Matrix elements for $<2,3,4>$ }\\ \hline

prime & ($v,x,y,w$) values \\
\hline
             103   &   (99,99,30,4), (43,43,48,60), (62,21,18,41),\ldots \\
             127   &  (65,65,19,62), (78,29,123,49), (53,117,27,74), \ldots \\
             151   & (133,79,15,18), (77,12,135,74), (52,15,21,99), \ldots \\
             263   & (82,82,11,181), (164,28,82,99), (191,16,54,72). \\
\hline
\end{tabular} \\

      The values ($v,x,y,w$) given are just a sample. There is a large number of quadruples that will work here. These were found by using the matrix representation for $<2,3,4> \times\,\, C_q$ given below and then raising the matrix $b$ to the power $q$.\\

            For the primes 103, 127, 151, and 263 (as well as for the other primes $p= 23$, 31,\ldots ) alternate matrix representations were found for the group $<2,3,4> \times\,\, C_q$. A selected sample of the values for ($x,y$) that appear in the matrix $b$ (which now has order $4*(p-1)/2)$,
\begin{equation}\label{A2:A2b.4}
b = \left(
       \begin{matrix}
           1 & x \\
           y & -1
        \end{matrix}
      \right),
\quad(b \text{ has order } 4*(p-1)/2),
\end{equation}
is given in the following table:\\

\begin{tabular}{|c|c|c|c|} \hline
\multicolumn{4}{|c|}{Table A6} \\ \hline
\multicolumn{4}{|c|}{ Matrix elements for $<2,3,4> \times \,\,C_q$ in (\ref{A2:A2b.4})}\\ \hline

prime & ($x,y$) & prime & ($x,y$) \\
\hline
              7 [4] &   (1,4)  &  79 [1] &  (15,22)  \\
                    &   (3,3)  &   103 * &  (1,44)   \\
                    &   (3,6)  &         & (100,44)  \\
                    &   (4,4)  &         &           \\
             23 [22]&   (1,9)  &   127 * &  (1,56)   \\
                    &   (1,16) &         &  (2,26)   \\
                    &   (2,7)  &         &  (24,2)   \\
             31 [16]&   (3,14) &   151 * &  (4,125)  \\
                    &   (4,26) &         &  (25,95)  \\
             47 [23]&   (1,19) &         &  (52,193) \\
                    &   (1,39) &         & (108,208) \\
                    &  (12,22) &  263 *  &   (1,29)  \\
                    &   (17,3) &         & (29,200) \\
             71 [8] &  (17,13) &         &           \\
                    &  (61,43) &         &           \\ \hline
\multicolumn{4}{|c|}{The numbers in [] brackets indicate the }\\
\multicolumn{4}{|c|}{number of solutions found. An $\ast$ means the} \\
\multicolumn{4}{|c|}{run was truncated before all solutions }\\
\multicolumn{4}{|c|}{were found. In case $p= 103$ all were }\\
\multicolumn{4}{|c|}{found but there was a very large number}\\
\multicolumn{4}{|c|}{and they were not counted. No computer}\\
\multicolumn{4}{|c|}{runs were done for the other primes. }\\ \hline
\end{tabular}\\

            The presentations for the groups $\textless2,3,4\textgreater  \times  C_q$
are given by (for $p= 7$):
\begin{equation}
\begin{split}
            a^3&=(a,b^2)=(a*b)^4*b^2=(a*(b^{-1}))^4*b^{-2} \\
               &=(a*b)^2*a^{-1}*b^{-1}*(a*b^{-1})^2*a^{-1}*b=1,
\end{split}
\end{equation}
and for $p > 7$,
\begin{equation}
\begin{split}
            a^3&=(a,b^2)=(a*(b^{-1}))^4*b^{-x} \\
               &=(a*b)^2*a^{-1}*b^{-1}*(a*(b^{-1}))^2*a^{-1}*b=1. \qquad \qquad \label{A2:A2b.3}
\end{split}
\end{equation}
Here $x = p-5$. These presentations for the group $\textless2,3,4\textgreater
\times C_q$ have been checked for primes up to 151 and for $p= 263$.\\

            This thus poses the interesting question in the theory
of group representations for the group $\textless2,3,4\textgreater$ :\\
\begin{quotation}
             What algebraic relationship does the pair ($x,y$) obey such that the presentation (\ref{A2:A2b.2}) (or (\ref{A2:A2b.3})) is satisfied? More generally, what conditions on the elements of the following order 4 matrix
\begin{equation}
b' = \left(
       \begin{matrix}
           v & x \\
           y & w
        \end{matrix}
      \right)
\end{equation}
are required in order for $<a,b'>$ to obey the presentation/rela-\linebreak tions (\ref{A2:A2b.2})? If one has access to a programming system such as Maple or Mathematica, then one may be able to investigate this problem rather easily; otherwise, it could be a fairly messy algebra problem. Remember what you are looking for is a representation valid for all primes $p\equiv  7$ mod(8).
\end{quotation}

\subsection{ Appendix 3. Matrix representations of the sylow 2-subgroups of $GL(n,p)$}

            The following describes how to obtain a matrix representation of the sylow 2-subgroups of the general linear groups, knowing what the corresponding sylow 2-subgroup is for $GL(2,p)$. The following is based upon the article by R. Carter and P. Fong \cite{8}.\\

            The sylow 2-subgroup of $GL(2,p)$, denoted by $S_2(2,p)$, is given by a pair of 2 by 2 matrices over the field $GF(p)$. Let these two matrices be denoted by $A$ and $B$. Then the sylow
2-subgroups of the higher-dimensional matrices are:\\

For the case of $GL(3,p)$:
\begin{align}
    &S_2(2,p) \times C_2 \quad \text{for $p \not\equiv 1$ mod(4),}\\
                     &S_2(2,p) \times C_4   \quad \text{for $p\equiv  1$ mod(4),}\\
                      & S_2(2,p) \times C_8 \quad \text{for $p\equiv  1$ mod(8),}\\
                       & S_2(2,p) \times C_{16} \quad \text{for $p\equiv  1$ mod(16)\ldots .}
\end{align}

      The matrix representations for these cases are just the direct products, i.e.,
\begin{equation}
M1 =
    \left(
      \begin{matrix}
           A & 0 \\
           0 & t
      \end{matrix}
  \right), \qquad
M2 = \left(
       \begin{matrix}
           B & 0 \\
           0 & t
        \end{matrix}
      \right),
\end{equation}
where $t$ is a generator of $C_2$, or $C_4$, or $C_8$  or $C_{16}$ ,\ldots , whichever
      is the appropriate choice depending upon $p$ of course.\\

            For the case of $GL(4,p)$ we have that the sylow 2-subgroups are wreath products of the corresponding $GL(2,p)$ sylow 2-subgroup with a $C_2$. The matrix representation for such a group is:
\begin{equation}
M1 =
    \left(
      \begin{matrix}
           A & 0 \\
           0 & I
      \end{matrix}
  \right), \qquad
M2 = \left(
       \begin{matrix}
           B & 0 \\
           0 & I
        \end{matrix}
      \right), \qquad
M3 = AI(4).
\end{equation} \qquad
      The matrix $AI$(4) is the following (anti-diagonal) $4 \times  4$ matrix:
\begin{equation}
   AI(4) = \left(
   \begin{matrix}
       0 & 0 & 0 & 1 \\
       0 & 0 & 1 & 0 \\
       0 & 1 & 0 & 0 \\
       1 & 0 & 0 & 0
    \end{matrix}
     \right),
\end{equation}
where $I$ is the $2 \times 2$ identity matrix in the block representation for the matrices $M1$ and $M2$.\\

            We shall continue this a little longer to show what the structure looks like for higher orders. For $GL(5,p)$ the sylow 2-subgroups are analogous to the case of $GL(3,p)$, namely the sylow
2-subgroups for $GL(4,p)$ cross a cyclic group of order $2^n$. For the sylow 2-subgroup of $GL(6,2)$,
$S_2(6,p)$, we have:
\begin{equation}
            S_2(6,p)  =  S_2(4,p) \times S_2(2,p),
\end{equation}
or
\begin{equation}
M1 =
    \left(
      \begin{matrix}
           A & 0 & | & 0 \\
           0 & I & | & 0 \\
           - & - & | & - \\
           0 & 0 & | & I
      \end{matrix}
  \right),\qquad
M2 = \left(
       \begin{matrix}
           B & 0 & | & 0 \\
           0 & I & | & 0 \\
           - & - & | & - \\
           0 & 0 & | & I
        \end{matrix}
      \right),
\end{equation}
\begin{equation}
M3 = \left(
      \begin{matrix}
                AI(4) & | & 0  \\
                 - & | & - \\
                0  & | & I
        \end{matrix}
      \right),
\end{equation}
\begin{equation}
M4 =
    \left(
      \begin{matrix}
             I(4) & | & 0 \\
             - & | & - \\
             0   & | & A
      \end{matrix}
  \right),\qquad
M5 = \left(
       \begin{matrix}
             I(4) & | & 0 \\
             - & | & - \\
              0   & | & B
        \end{matrix}
      \right),
\end{equation}
where $I(4)$ is the $4 \times 4$ identity matrix. This representation can probably be simplified to fewer generators, but this gives the overall structure of the groups.\\

            For the case of $GL(7,p)$ we have just the case of $GL(6,p)$ cross the appropriate cyclic group again. For the case of $GL(8,p)$ we have a second wreath product occurring, namely [$S_2(2,p)wrC_2 ] wrC_2$ . The matrix form here is:

\begin{equation}
M1 =
    \left(
      \begin{matrix}
           A & 0 & | & 0 \\
           0 & I & | & 0 \\
           - & - & | & - \\
           0 & 0 & | & I(4)
      \end{matrix}
  \right),\qquad
M2 = \left(
       \begin{matrix}
           B & 0 & | & 0 \\
           0 & I & | & 0 \\
           - & - & | & - \\
           0 & 0 & | & I(4)
        \end{matrix}
      \right),
\end{equation}
\begin{equation}
M3 =
    \left(
      \begin{matrix}
              AI(4)& | & 0 \\
              - & | & - \\
              0 & |& I(4)
      \end{matrix}
  \right),\qquad
M4 =AI(8),
\end{equation}
      where $AI(8)$ is the $8 \times  8$ anti-diagonal matrix:
\begin{equation}
M4 =
    \left(
      \begin{matrix}
           0 & 0 & 0 & 0 & 0 & 0 & 0 & 1 \\
           0 & 0 & 0 & 0 & 0 & 0 & 1 & 0 \\
           0 & 0 & 0 & 0 & 0 & 1 & 0 & 0 \\
           0 & 0 & 0 & 0 & 1 & 0 & 0 & 0 \\
           0 & 0 & 0 & 1 & 0 & 0 & 0 & 0 \\
           0 & 0 & 1 & 0 & 0 & 0 & 0 & 0 \\
           0 & 1 & 0 & 0 & 0 & 0 & 0 & 0 \\
           1 & 0 & 0 & 0 & 0 & 0 & 0 & 0
      \end{matrix}
  \right).
\end{equation}

      For higher orders the sequence continues in this manner, i.e., whenever we have $GL(2^n,p)$ we have an additional wreath product, which is represented by a  $2^n  \times 2^n$  anti-diagonal matrix, acting on the previous wreath product.\\

Other cases proceed in like manner, e.g., $GL(15,p)$'s sylow 2-subgroup will have the schematic block structure:
\begin{equation}
M =
    \left(
      \begin{matrix}
           S_2(8,p) & | & 0 & | & 0 &| & 0 \\
            - & | & - & | & - & | & - \\
            0 & | & S_2(4,p) & | & 0 & | & 0  \\
            - & | & - & | & - & | & - \\
            0 & | & 0 & | & S_2(2,p) & | &0 \\
            - & | & - & | & - & | & - \\
            0 & | & 0 & | & 0 & | & t
      \end{matrix}
  \right),
\end{equation}
where $t$ is a generator of the appropriate cyclic group of order
$2^n$ (see the $GL(3,p)$ case above).\\

\section{                     Acknowledgements}

This work could not have been done without the aid of a number of individuals associated with Michigan State University (in the early 1980s), the University of Rhode Islands' Engineering Computer Laboratory (1984), the High Energy Particle Physics group at Syracuse University (from 1985 to about 1990) and Brown University's Cognitive and Linguistic Sciences Department computers (from 1990 on). The long list of acknowledgements given in the previous report on the groups of orders $8p$ and $8p^{2}$ also apply to this work as well. The major difficulty the authors have had over most of this span of time is the lack of personal contact with group theorists and others involved with computational group theory which would have made the effort more interesting as well as improving the presentations given in these papers. The authors would appreciate any comments or suggestions on the mode of presentation or requests for inclusion of additional information on these groups or their automorphism groups.\\

Dr. M. F. Fry helped us in the early stages of this work to obtain a complete list of the groups of order 144.\\

            The early version of CAYLEY we were using at Michigan State
      University would not accept presentations, but would accept a
      permutation input for the groups. Dr. Richard Hartung (then with
      the Physics Department) allowed us to use his Todd-Coxeter program
      not only to check the relations for our groups but also
      modified his program so that we could get permutation representations
      for the relations that we were using. Without his assistance much
      of the early work done at Michigan State University would not have
      been possible.\\

            The work on the higher-order cases was done on the Brown
      University computers, but most of the other cases were done at the
      other locations.\\

\begin{tabular}{|c|c|c|} \hline
  \multicolumn{3}{|c|}{Table 1} \\ \hline
  \multicolumn{3}{|c|}{Groups of order 16}  \\ \hline
extension & relations & automorphism group \\ \hline
$C_{16}$ & $a^{16}=1$ & $C_4 \times C_2 $\\ \hline
$C_8 \times  C_2$ & $ a^8=b^2=(a,b)=1$ & $ D_4 \times  C_2$\\ \hline
$C_4 \times  C_4$ &  $a^4=b^4=\left(a,b\right)=1$ & group of order 96, \\
& & degree 8, group \\ & & generators: \\
 & &  $a=$(1,3,2)(5,6)(7,8), \\ & &
$b=$(1,2,3,4)(5,6,7,8),\\
   & &   $A_4 @ D_4, C_4$ acts as $C_2$ \\ \hline
$C_4 \times  C_2 \times  C_2$& $ a^4=b^2=c^2=(a,b)=(a,c)$ & group of order
192\\
   & $ =(b,c)=1$   &  $ \left( \left( Q_2\text{Y}Q_2\right) @ C_3 \right) @ C_2
$ \\ \hline
$C_2 \times  C_2 \times  C_2 \times  C_2$ & $ a^2=b^2=c^2=d^2=(a,b)=$&
$GL(4,2) = A_8$ \\
  &      $(a,c)=(a,d)=(b,c)=(b,d)$ & \\
   & $=(c,d)=1$ & \\ \hline
   $      D_4\times  C_2 $&$  a^4=b^2=a^b*a= $  &   Hol($C_4 \times  C_2$)   \\
                     &$ c^2=(a,c)=(b,c)=1$      &    (\#259 order 64)    \\ \hline
$     Q_2 \times  C_2   $&$ a^4=b^4=a^2*b^2=a^b*a=$       &   order 192
  \\
                     &$ c^2=(a,c)=(b,c)=1  $    &$   ((C_4 \times  C_4) @ C_3)
$\\
                      &                        &$      @ (C_2 \times  C_2)
  $ \\ \hline
   $      C_4$Y$Q_2 $ &$ a^2=b^2=c^4=a^b*c^2*a=    $   &$   S_4 \times  C_2
   $ \\
                     &$    a^c*a=b^c*b=1        $&                          \\
\hline
      $(4,4 |2,2)$    &$ a^4=b^2=c^2=(a,b)=(b,c)  $ &$   (C_2 \times  C_2) wr
C_2      $\\
   $ (C_4 \times  C_2) @ C_2 $&$  =a^c*a^{-1}*b=1 $  &      (\#33 order 32)
\\ \hline
      $<2,2 | 4;2>$  &$ a^4=b^4=a^b*a=1 $  &$   (C_2 \times  C_2) wr C_2
$\\
   $  C_4 @ C_4      $&                          &        (\#33 order 32)  \\ \hline
         $<2,2|2>$ &$ a^8=b^2=a^b*a^{-5}=1   $    &$    D_4 \times  C_2
$\\
   $  C_8 @ C_2      $&                         &                       \\
\hline
   $  D_8            $&$ a^8=b^2=a^b*a=1 $         &    Hol($C_8$)
\\
   $  C_8 @ C_2        $&                          &
\\ \hline
      $<-2,4|2>$       &$ a^8=b^2=a^b*a^{-3}=1  $     &$    D_4 \times  C_2
  $ \\
      $C_8 @ C_2        $ &                          &                      \\
\hline
$     Q_4      $      &$ a^8=b^4=a^4*b^{-2}=       $  &   Hol($C_8$)            \\
        $<2,2,4>$        &$       a^b*a=1           $&                         \\
\hline
\end{tabular}\\

\begin{tabular}{|c|}\hline
Notes for Table 1\\ \hline
     relations for Aut($C_4 \times  C_4)$: $A_4 @ D_4$ version.  \\
             $ a^2=b^3=b*a^3=c^4=d^2=c^d*c=(a,d)=(b,d)=$                 \\
                 $(a,c)=b^c*b=1 $;                                       \\
                (note not same as given by above permutations!)        \\
            relations of above permutations:                             \\
                   $ a^6=b^4={a*b}^2=(a*b^{-1})^4=$                        \\
           $ a^3*b*a^{-1}*b^{-1}*a*b^{-1}*a^{-1}*b=1$;           \\ \hline
relations for Aut($C_4 \times  C_2 \times  C_2$):   \\
            $ b^4=c^2*a*c^{-1}*a=c*a^2*c*a^{-1}=b^2*a^{-3}=$       \\
                 $ c*b*c*b^{-1}*c*b^{-1}*a=d^2=a^d*((b*a*c)^{-1})=$  \\
                $ b^d*a^{-1}*c=c^d*((a*b*c)^{-1})=1$;                 \\ \hline
            relations for Aut($Q_2 \times  C_2$):                        \\
              $a^4=b^4=(a,b)=c^3=a^c*b^{-1}*a^{-1}=b^c*b^2*a^{-1}=$ \\
              $     d^2=e^2=(d,e)=                               $ \\
               $    a^d*b*a^2=b^d*a*b^2=c^d*c*a^2=              $ \\
             $ a^e*b^{-1}*a^2=b^e*b^2*a^{-1}=c^e*c=1.    $ \\   \hline
               \\
            Also note that relations for the automorphism groups with  \\
            orders 32 or 64 can be found in Sag and Wamsley's article  \\
            \cite{4} on the minimal presentations of the groups of order    \\
            $2^n$ with $n$ less than or equal to 6. The presentations used  \\
            in \cite{1} were taken from this article.                      \\
             \\  \hline
\end{tabular}

\begin{tabular}{|c|c|c|c|c|c|c|} \hline
  \multicolumn{7}{|c|}{Table 2a} \\ \hline
  \multicolumn{7}{|c|}{Groups of order $16p$ and their automorphism groups}
\\ \hline
2-group&$C_2$&\multicolumn{4}{|c|}{image of 2-group}&$p=3$, 5, or 7* \\
\cline{3-6}
    & action&$C_2$&$C_4$&$C_8$&$C_{16}$& (normal 2-group) \\ \hline
$  C_{16}  $&  $a$   &$C_4 \times  C_2  $&$[a]   C_4$&$ C_2$& I    &    \\
\hline
$C_8 \times  C_2$& $a$ &$D_4 \times  C_2 $ & $[a]  D_4$&$ C_2 $&  & \\
\cline{2-7}
         &  $b$    & $  1^3$&  $[a,b]  D_4 $&   &   &           \\ \hline
$C_4 \times  C_4$ &  $a$ & $1^2 wr C_2$& $ [a] D_4 $& &  & order 384     \\
\hline
$C_4 \times C_2$&  $a$    &Aut($2,1^2$)& $[a]  S_4$&   & $  S_4 \times
C_2 $&$ S_4 \times C_2 $ \\ \cline{2-7}
$\times  C_2   $&  $b$    &$1^2 wr C_2 $&   &  & &         \\ \hline
$1^4 $&  $a$    & Hol($1^3$) &       &  &  &$S_4\times S_3$ \\\cline{2-7}
      &    &     &       &  & &order (5760)  \\ \cline{2-7}
       &     &       &  & & &$p=5$:(complete)  \\ \cline{2-7}
      &     &     &  &  & & $p=7$: (complete  \\
      &     &     &  &  & &group of order 168)  \\ \hline
$ D_4 \times  C_2$&  $b$    & Hol($2,1$) &           & & &
\\ \cline{2-7}
         &  $a$    &$1^2 wr C_2 $&     &  &  &           \\ \cline{2-7}
         &  $c$   &$D_4 \times  C_2  $&           &  &  &              \\
\hline
$Q_2\times C_2 $&  $a$ &$ S_4 \times  C_2   $&    & &  &$ S_4 \times  C_2 $
\\ \cline{2-7}
         &  $b$    & Hol($2,1$) &       &  & &        \\ \hline
$C_4$Y$Q_2 $& $a,b,c$ &$ S_4 \times  C_2 $&      &  &  &$ S_4 \times  C_2  $
\\ \cline{2-7}
         &  $b,c$  &$ D_4 \times  C_2  $&           & & &  \\ \cline{2-7}
         &  $a,c$  &$ D_4 \times  C_2 $&        &  &  &       \\ \hline
$(4,4|2,2)$&  $b$   &$1^2 wr C_2 $&$ [a]    D_4 $&  & &          \\ \cline{2-7}
         &   $c$   &$   1^4    $&         &     &     &        \\ \hline
$<2,2|4,2>$&  $b$   &$1^2 wr C_2 $&$ [b]    D_4 $& &  &           \\ \cline{2-7}
         &   $a$   &$   1^4  $  &           & & &   \\ \hline
$<2,2|2>$ &   $a$   &$ 1^2 wr C_2$&$ [a]    D_4 $&  &  &            \\
\cline{2-7}
       &   $b$   &  $ 1^3 $&   $ [a,b]  D_4 $& &  &             \\ \hline
$    D_8  $&   $b$   & Hol($C_8$) &       &  &  &             \\ \cline{2-7}
         &   $a$   &$ D_4 \times  C_2  $&     & &&     \\ \hline
$   QD_8  $&   $a$  &$ D_4 \times  C_2 $&   &  &  &           \\ \cline{2-7}
    or   &   $b$   &$ D_4 \times  C_2  $&   &  & &           \\ \cline{2-7}
$<-2,4|2>$&  $a,b$  &$ D_4 \times  C_2  $&   && &    \\ \hline
$   Q_4 $&  $b$     & Hol($C_8$) &      &  &  &            \\ \cline{2-7}
         &  $a$  &$ D_4 \times  C_2  $& && &       \\ \hline
   Totals &            & 28 &         9 &        2&     1&   6 + 4(in Table 2b)
\\ \hline
\multicolumn{7}{|c|}{* See Tables 2b and 2c for more details on these
cases.}\\ \hline
\end{tabular}

\begin{tabular}{|l|} \hline
Notes to Table 2a \\ \hline
            In the $C_4$ image column the [--] gives the 2-group \\
                generators acting on the $C_p$-subgroup, and the \\
                other entry (e.g., $C_4$) gives the automorphism \\
                groups' invariant factor as in the previous tables. \\ \hline
            $1^2$ means $C_2 \times  C_2.$ \\
            $1^3$ means $C_2 \times  C_2 \times C_2.$ \\
            $1^4$ means $C_2 \times  C_2\times C_2 \times  C_2.$  \\
            $1^2 wr C_2$ means $(C_2 \times C_2) wr C_2 $ (wreath product).\\
            Aut(2,$1^2$) means automorphism group of $(C_4 \times C_2 \times  C_2).$ \\
            Hol(2,1) means the holomorph of $(C_4 \times  C_2).$\\
            Hol$(C_8)$  means the holomorph of $C_8 $.\\ \hline
      Note for $p\equiv  1$ mod(2) only, we just have the $C_2$ image cases.\\
      If $p\equiv  1$ mod(4), then we have the $C_2$ and $C_4$ image cases.\\
      If $p\equiv  1$ mod(8), then we have the $C_2$, $C_4$ and $C_8$ image cases.\\
      If $p\equiv  1$ mod(16)....etc., then we have the $C_2$, $C_4$, $C_8$,\\
      and $C_{16}$ image cases.\\
\hline
\end{tabular}

\begin{tabular}{|c|c|c|} \hline
  \multicolumn{3}{|c|}{Table 2b} \\ \hline
  \multicolumn{3}{|c|}{For the prime $p=3$ the following special groups arise.}
  \\ \hline
group & presentation & automorphism group \\ \hline
       $  (C_4 \times  C_4) @ C_3 $   &$   a^4=b^4=(a,b)=c^3=   $& complete
group     \\
                             &$ a^c*b^{-1}=b^c*a*b=1    $&   of order 384
  \\ \hline
       $[(C_2 \times  C_2) @ C_3] \times  C_4  $&        &   $      S_4
\times  C_2 $   \\
        or\quad$ A_4 \times  C_4           $&            &
\\ \hline
       $   A_4 \times  C_2 \times  C_2       $&     & $  S_4 \times  S_3 $
\\ \hline
       $ (C_2 \times  C_2 \times  C_2 \times  C_2)@C_3$&$
a^2=b^2=c^2=d^2=(a,b)$&   complete group  \\
                            &$ =(a,c)=(a,d)=(b,c)   $&     of order 5760    \\
                              &$ =(b,d)=(c,d)=e^3 $&                    \\
                             &$ =a^e*b=b^e*a*b $&              \\
                              &$ =c^e*d=d^e*c*d=1 $&                  \\
\hline
       $ (Q_2 @ C_3) \times  C_2$&    &            $   S_4 \times  C_2 $  \\
       or $SL(2,3) \times  C_2      $ &                      &
\\ \hline
       $ (C_4$Y$Q_2) @ C_3   $&$ a^2=b^2=c^4=(a,b)*c^2=   $&$       S_4
\times  C_2   $ \\
                          &$ (a,c)=(b,c)=d^3=a^d*b  $&
\\
                           &$ =b^d*c^{-1}*b*a=(c,d)    $&                  \\
                            & $ =(a*b*c)^d*a=1          $ &            \\
\hline
       $  SL(2,3) @ C_2   $&                   &    $       S_4 \times  C_2
$     \\
          or $GL(2,3)$         &                       &                  \\
\hline
         $<2,3,4>$            &$ a^2=b^3=c^4=a*b*c    $&$     S_4 \times  C_2  $
\\ \hline
       $    S_4 \times  C_2   $ &    &                $       S_4 \times  C_2
   $\\ \hline
       $    A_4 @ C_4   $&$ a^2=b^3=(b*a)^3=c^4=   $&$       S_4 \times  C_2
$   \\
                         &    $ a^c*b^{-1}*a^{-1}*b=   $ &
\\
                          &   $ b^c*b*a^{-1}=1           $&             \\
\hline
\multicolumn{3}{|c|}{The first six cases above have a normal sylow 2-subgroup.} \\
\multicolumn{3}{|c|}{The last four cases do not have a normal sylow
subgroup.}\\ \hline
\multicolumn{3}{|c|}{Presentation for order 384 automorphism group:  } \\
\multicolumn{3}{|c|}{$ a^4=b^3=c^4=a^2*b*c^{-2}*b^{-1}=a*b*a^{-1}*c*b*c=
$}\\
\multicolumn{3}{|c|}{$ a*b^{-1}*a*c*b^{-1}*c=a*c^{-1}*b*c*a^{-1}*b=     $
}\\
\multicolumn{3}{|c|}{$
(a*b)^2*(a^{-1}*b)^2=a*c*a*c^{-1}*a*c^{-1}*a^{-1}*c=1$ }\\ \hline
\multicolumn{3}{|c|}{Presentation for the order 5760 automorphism group: }
\\
\multicolumn{3}{|c|}{$ a^2=(a*b)^2=b^4=c^6=(a*c)^4=a*c^3*a*c^{-3} =  $}\\
\multicolumn{3}{|c|}{$ b*c^3*b*c^{-3}=a*b*c*b*c^2*b^{-1}*c*a*c= $} \\
\multicolumn{3}{|c|}{$a*b^2*c^{-1}*b^2*c*a*c^{-1}*b^{-2}*c= $}  \\
\multicolumn{3}{|c|}{$a*b^2*c*b^2*c^2*a*c*b^{-2}*c^{-1}=1 $}\\ \hline
\end{tabular}

\begin{tabular}{|c|c|c|} \hline
  \multicolumn{3}{|c|}{Table 2c} \\ \hline
  \multicolumn{3}{|c|}{For the primes 5 and 7 we have one case each in which
we have
}  \\
\multicolumn{3}{|c|}{a normal sylow 2-subgroup. These cases are listed
below. }\\ \hline
               group &       presentation  &  automorphism   \\
   & & group  \\ \hline
      $(C_2 \times  C_2 \times  C_2 \times  C_2)@C_5$&$
a^2=b^2=c^2=d^2=(a,b)= $&   complete group   \\
                           &$ (a,c)=(a,d)=(b,c)=    $&     of order 960
\\
                            &$ (b,d)=(c,d)=e^5=       $&                      \\
                            &$  a^e*a*c*d=b^e*a=       $&                     \\
                            &$ c^e*a*b*c*d=   $&                      \\
                            &$ d^e*a*d=1               $&                  \\
                             &                       &                      \\
                             & or                    &                    \\
                             &$   a^5=b^2=c^5=c*a*b=1   $&
\\
                            &                       &                 \\
\hline
$ (C_2 \times  C_2 \times  C_2) @ C_7 $&$ a^2=b^2=c^2=d^2=(a,b)$&  complete
group of \\
      $      \times  C_2            $&$ =(a,c)=(a,d)=(b,c)=  $&       order 168
    \\
                             &$ (b,d)=(c,d)=e^7=      $&                      \\
                             &$ a^d*b=b^d*c=c^d*a*b  $&                      \\
                             &$     =(d,e)=1      $&                    \\ \hline
\multicolumn{3}{|c|}{The order 168 group is the complete group of this
order.} \\
\multicolumn{3}{|c|}{This group is isomorphic to a degree 8 permutation
group } \\
\multicolumn{3}{|c|}{with generators:  } \\
\multicolumn{3}{|c|}{$a$ =  (1,2)(3,6,7,4,5,8),   $b$ =  (1,7,2,6,4,5)(3,8)
} \\
\multicolumn{3}{|c|}{and presentation:                           } \\
\multicolumn{3}{|c|}{$  a^6=(a*b^{-1})^3=b^6=a^2*b^{-1}*a*b*a^{-2}*b= $
} \\
\multicolumn{3}{|c|}{$  (a^2*b)^2*a^{-1}*b^{-2}=1. $               } \\
\hline
\multicolumn{3}{|c|}{A presentation for the complete group of order 960 is
} \\
\multicolumn{3}{|c|}{$  b^2=a^4=a*c^2*a^{-1}*c^{-1}=(a,b)=(a^2*c^{-3})^2=  $
} \\
\multicolumn{3}{|c|}{$      c*b*c^3*b*c^{-4}*b=1. $               } \\ \hline
\end{tabular}

\begin{tabular}{|c|c|c|c|c|} \hline
  \multicolumn{5}{|c|}{Table 3a} \\ \hline
  \multicolumn{5}{|c|}{Basic groups of order $16p^2$ and their automorphism
groups*}  \\ \hline
2-group & $C_2\times C_2$& \multicolumn{2}{|c|}{image of 2-group}&  \\
\cline{3-4}
   & $p$-factor &$C_2\times C_2 $&$ C_4$& order 8 and 16  \\
  & (if any)& & & cases (see notes  \\
   &         & & & for details) \\ \hline
$  C_{16} $&              &none here &$ a: C_4 \times  [144] $&$ C_8:  \text{order }
288  $\\ \hline
$C_8 \times  C_2$ &$S_3 \times  S_3 (a,b)$&$  1^3$&$ a: D_4 \times  [144]$
&$ C_8:        $\\ \cline{2-4}
    &     $(ab,b)$   &[576](54,4) &$ ab:D_4 \times  [144] $&$ C_2 \times
[144] $  \\ \hline
$C_4 \times  C_4 $&  $(a,b)$&$ (S_3 \times  1^2) $&$ a: G=[36] \times  C_4$&
none      \\
          &              &  $  wr C_2 $ &Aut $ = D_4 \times  [144]$& \\
\hline
$C_4 \times  C_2  $&    $(b,c)$     &[576](54,4) &$ a: [36] \times  1^2 $&
none  \\ \cline{2-3}
  $\times  C_2 $& $S_3 \times  S_3 (a,b) $&   \#33     &  Aut = &         \\
\cline{2-3}
          &   $(ab,b)$   &[2304](90,2)& $   S_4 \times  [144] $ &\\ \hline
$ (1^4)   $&    $(a,b)$     & [6912]     &  no $C_4$ cases & none    \\ \hline
$D_4 \times  C_2 $& $S_3 \times  S_3 (a,c) $&$    1^3    $&    & for $p= 3, G =$
  \\ \cline{2-3}
   &$S_3 \times  S_3 (b,c) $&$ D_4 \times  C_2  $&    no &$(S_3wrC_2) \times  C_2 $ \\ \cline{2-3}
        &$S_3 \times  S_3 (ab,b)$&$   \#33  $&$    C_4   $& ($D_4$ image)
and    \\ \cline{2-3}
        &  $(abc,bc)$   &[576](54,4) &   cases       & Aut =
\\ \cline{2-3}
         &  $(abc,ac)$   &[1152](81,2) &              & order 576(27,2)
\\ \cline{2-3}
         & $(ab,a)$     &[2304](90,2) &               &
\\ \hline
$Q_2 \times  C_2$&$S_3 \times  S_3 (b,c) $&$   D_4 \times  C_2  $&    no
&$ Q_2  $ \\ \cline{2-3}
         &  $(c,bc)$     &[1152](81,2) &$    C_4     $&               \\
\cline{2-3}
         &  $(ab,a)$     &[2304](90,2) &   cases      & Aut =
\\
         &               &             &              & order 1728(20,1) \\ \hline
$C_4$Y$Q_2 $&$S_3 \times  S_3 (ac,a)$&$   1^3   $&        & none here \\
\cline{2-3}
       &$S_3 \times  S_3 (x,ab)$&$   D_4 \times  C_2  $&     no        &  \\
\cline{2-3}
        &$S_3 \times  S_3 (x,bc)$&$   D_4 \times  C_2 $&$     C_4   $ &\\
\cline{2-3}
       &$S_3 \times  S_3 (b,ac)$&$   D_4 \times  C_2  $&    cases   &   \\
\cline{2-3}
       &  $(ab,a)$     &[576](54,4) &               &             \\
\cline{2-3}
        &  $(bc,c)$      &[576](54,4) &               &             \\ \hline
$(4,4|2,2)$&$S_3 \times  S_3 (a,b) $&$    1^4    $& $b$:     &$ D_4     $\\
\cline{2-3}
          & $(b,ab)$     &$(S_3 \times  1^2) $&$ D_4 \times [144] $&  Aut =
\\
          &              &$   wr C_2   $&            &$   S_3 wr C_2 \times
1^2 $ \\ \hline
$<2,2|4,2>$&$S_3 \times  S_3 (ab,b)$&$   1^4  $& $b$:   &$ D_4$: Aut =      \\
          &              &            &               &   order 576(27,2)
\\ \cline{2-3}
          &  $(ab,a)$      &$ (S_3 \times  1^2)$&$  D_4 \times  [144] $&$
Q_2$: Aut =         \\
          &              &$   wr C_2    $&              &     order
576(24,2)\\ \hline
$<2,2|2>$  & $(ab,b)$      &[576](54,4) & 2 both        & none here
\\ \cline{2-3}
       &$S_3 \times  S_3 (ab,a)$&$    1^3   $&$ D_4 \times  [144]  $  &   \\
        &              &            & $[a]$,$[ab]$      &               \\
\hline
\multicolumn{5}{|c|}{Table 3a continued on next page} \\ \hline
\end{tabular}

\begin{tabular}{|c|c|c|c|c|} \hline
  \multicolumn{5}{|c|}{Table 3a \quad continued} \\ \hline
  \multicolumn{5}{|c|}{Basic groups of order $16p^2$ and their automorphism
groups*}  \\ \hline
2-group & $C_2\times C_2$& \multicolumn{2}{|c|}{image of 2-group}&  \\
\cline{3-4}
   & $p$-factor &$C_2\times C_2 $&$ C_4$& order 8 and 16 cases \\
  & (if any)& & & (see notes for details) \\ \hline

  $ D_8 $&$S_3 \times  S_3 (a,b)$&$   D_4 \times  C_2  $&$  \text{ no }  C_4
$&$ D_4 $  \\ \cline{2-3}
          &  $(ab,a)$      &[1152](63,2)&   cases       & Aut order =
576(27,2) \\ \hline
$ C_8 @ C_2$&$S_3 \times  S_3 (ab,b) $&$ D_4 \times  C_2  $&        &$D_4:
$     \\ \cline{2-3}
           &           &    &     no        & Aut = $S_3 wr C_2 \times
1^2$\\
$QD_8$          &              &           &$     C_4        $&
  \\ \
  &$S_3 \times  S_3 (ab,a)$&$ D_4 \times  C_2   $&    cases   &$C_8 @ C_2
---> $ \\ \cline{2-3}
$<-2,4|2>$&$ S_3 \times  S_3 (a,b) $&$ D_4 \times  C_2    $&               &
  complete group.   \\ \hline
$  Q_4 $&$S_3 \times  S_3 (a,b) $&$  D_4 \times  C_2  $&  no $C_4$
&$D_4: $ Aut =    \\ \cline{2-3}
       &  $(ab,a)$      &[1152](63,2)&  cases        &   order 576(27,2)  \\
\hline
  Totals  &            &  35      &    9        &$  C_8=2, D_4=6, Q_2=2,
$\\
           &              &            &               & $QD_8=1$ for $p\equiv 3$ mod(8) \\
\hline
\multicolumn{5}{|c|}{}\\
\multicolumn{5}{|c|}{* For $p$ greater than 3, $S_3$ goes over into
Hol($C_p$). } \\
\multicolumn{5}{|c|}{$x$ in the $C_2 \times  C_2$ $p$-factor column stands for
$abc$, e.g., in $C_4 \text{Y} Q_2$,     } \\
\multicolumn{5}{|c|}{($x,ab$) = ($abc,ab$).  } \\
\multicolumn{5}{|c|}{}\\ \hline
\multicolumn{5}{|c|}{}\\
\multicolumn{5}{|c|}{ The other cases, e.g., (1152)(63,2) are discussed
individually } \\
\multicolumn{5}{|c|}{in the text, primarily in section 3.3.} \\
\multicolumn{5}{|c|}{}\\ \hline
\multicolumn{5}{|c|}{}\\
\multicolumn{5}{|c|}{[144] means the complete group $(C_3 \times C_3)@(C_8 @ C_2)$ given by}\\
\multicolumn{5}{|c|}{equation (\ref{CG144}) above. Its presentation is:}\\
\multicolumn{5}{|c|}{$a^8=b^2=a^b*a^5=p^3=q^3=(p,q)$}\\
\multicolumn{5}{|c|}{$=p^a*q*p^{-1}=q^a*q^{-1}*p^{-1}=(p,b)=q^b*q=1$}\\
\multicolumn{5}{|c|}{or, in terms of permutations,}\\
\multicolumn{5}{|c|}{$a = (2,9,4,6,3,7,5,8),\quad   b = (4,5)(6,9)(7,8),$}\\
\multicolumn{5}{|c|}{$p = (1,3,2)(4,7,6)(5,8,9),\quad   q = (1,5,4)(2,9,6)(3,8,7).$} \\
\multicolumn{5}{|c|}{}\\ \hline
\multicolumn{5}{|c|}{}\\
\multicolumn{5}{|c|}{The last column gives the group with a $D_4, C_8, Q_2,$ or $QD_8$}\\
\multicolumn{5}{|c|}{action on the $C_p \times C_p$ group which is the normal subgroup here.}\\
\multicolumn{5}{|c|}{}\\ \hline
\end{tabular}

\begin{tabular}{|c|c|c|c|c|c|c|} \hline
\multicolumn{7}{|c|}{Table 3b} \\ \hline
\multicolumn{7}{|c|}{Automorphism groups of the groups of order $16p^2$} \\
\multicolumn{7}{|c|}{arising from a $D_4$ or a $Q_2$ action} \\ \hline
2-group & image & \multicolumn{5}{|c|}{factor acting on the $p$-group $ C_p
\times C_p $ } \\ \cline{3-7}
           &        &   $p=3$    &   $p=5$    &    $p=7$  &   $p=11$   &   $p=13$  \\ \hline
$D_4 \times  C_2$ &$   D_4  $& $<\!-2,4|2\!>$ &$ C_4 wr C_2$ &$ D_8 \times
C_3$ & $<\!-2,4|2\!>$ &$ C_4 wr C_2$\\
& & & & &$\times C_5$ &$ \times C_3$\\ \hline
$C_4 @  C_4$ &$   D_4  $& $<\!-2,4|2\!>$ &$ C_4 wr C_2$ &$ D_8 \times  C_3$ &
$<\!-2,4|2\!>$ &$ C_4 wr C_2$\\
& & & & &$\times C_5$ &$ \times C_3$\\ \hline
$C_4 @ C_4$ &$   Q_2  $& $<\!-2,4|2\!>$ &$ C_4 wr C_2$ &$ Q_4 \times
C_3$ & $<\!-2,4|2\!>$ &$ C_4 wr C_2$\\
& & & & &$\times C_5$ &$ \times C_3$\\ \hline
$D_8$ &$   D_4  $& $<\!-2,4|2\!>$ &$ C_4 wr C_2$ &$ D_8 \times  C_3$ &
$<\!-2,4|2\!>$ &$ C_4 wr C_2$\\
& & & & &$\times C_5$ &$ \times C_3$\\ \hline
$Q_4 $ &$   D_4  $& $<\!-2,4|2\!>$ &$ C_4 wr C_2$ &$ D_8 \times  C_3$ &
$<\!-2,4|2\!>$ &$ C_4 wr C_2$\\
& & & & &$\times C_5$ &$ \times C_3$\\ \hline
\multicolumn{2}{|c|}{order of Aut($g$) =} &576 &3200 & 9408 &   38720  & 64896
\\ \hline
\multicolumn{7}{|c|}{ } \\ \hline
2-group & image & \multicolumn{5}{|c|}{factor acting on the $p$-group $ C_p
\times C_p $ } \\ \cline{3-7}
           &        &   $p=17$    &   $p=19$    &    $p=23$  &   $p=29$   &   $p=31$  \\
\hline
$D_4 \times  C_2$ &$ D_4$ &order 128 & $<\!-2,4|2\!>$ &$ D_8 \times  C_{11}$ &$
C_4 wr C_2$ &   $D_8$   \\
         &        & ncl = 56 & $ \times  C_9 $  &      & $ \times  C_7$ &$
\times  C_{15}$ \\ \hline
$C_4 @  C_4$ &$ D_4$ &order 128 & $<\!-2,4|2\!>$ &$ D_8 \times  C_{11}$ &$ C_4
wr C_2$ &   $D_8$   \\
         &        & ncl = 56 & $ \times  C_9 $  &      & $ \times  C_7$ &$
\times  C_{15}$ \\ \hline
$C_4 @  C_4$ &$ Q_2$ &order 128 & $<\!-2,4|2\!>$ &$ Q_4 \times  C_{11}$ & &
\\
         &        & ncl = 56 & $ \times  C_9 $  &      &  & \\ \hline
$D_8 $ &$ D_4$ &order 128 & $<\!-2,4|2\!>$ &$ D_8 \times  C_{11}$ &$ C_4 wr C_2$
&   $D_8$   \\
         &        & ncl = 56 & $ \times  C_9 $  &      & $ \times  C_7$ &$
\times  C_{15}$ \\ \hline
$Q_4 $ &$ D_4$ &order 128 & $<\!-2,4|2\!>$ &$ D_8 \times  C_{11}$ &$ C_4 wr C_2$
&   $D_8$   \\
         &        & ncl = 56 & $ \times  C_9 $  &      & $ \times  C_7$ &$
\times  C_{15}$ \\ \hline
\multicolumn{2}{|c|}{order of Aut($g$) =} &147968 & 207936 & 372416  & 753536
& 922560   \\ \hline
\multicolumn{7}{|c|}{The general structure of these automorphism groups
appears to be } \\
\multicolumn{7}{|c|}{} \\
\multicolumn{7}{|c|}{$U \cong (C_p \times  C_p \times  C_2 \times  C_2) @ X $}\\
\multicolumn{7}{|c|}{}\\
\multicolumn{7}{|c|}{ where $X$ is the factor in the above table and $U$ has order
$2^5p^2(p-1)$. }\\
\multicolumn{7}{|c|}{The abbreviation ncl in the above table means the number }\\
\multicolumn{7}{|c|}{of conjugacy classes in the group of order 128. This group }\\
\multicolumn{7}{|c|}{is \#902 in the small group library and has a presentation }\\
\multicolumn{7}{|c|}{given by $a^2=a*c*a*c^{-1}=b^4=(b,c)=(a*b)^2(a*b^{-1})^2$}\\
\multicolumn{7}{|c|}{$=(a*b)^2*c^{-4}=1.$}\\
\multicolumn{7}{|c|}{}\\
\multicolumn{7}{|c|}{A direct determination of the $p\equiv  9$ mod(16) cases was
not possible}\\
\multicolumn{7}{|c|}{at this site with CAYLEY. The orders of these groups were
too large}\\
\multicolumn{7}{|c|}{to handle here. An explicit presentation for most of these
groups}\\
\multicolumn{7}{|c|}{can be found in section 3.3 of the text. }\\ \hline
\end{tabular}

\begin{tabular}{|c|c|} \hline
\multicolumn{2}{|c|}{Table 4a} \\
\multicolumn{2}{|c|}{ $C_4$ actions} \\ \hline
\multicolumn{2}{|c|}{} \\
\multicolumn{2}{|c|}{If the automorphism group in the} \\
\multicolumn{2}{|c|}{$p=3$ case is $D_4 \times$ [144],} \\
\multicolumn{2}{|c|}{then the automorphism groups in the} \\
\multicolumn{2}{|c|}{$p\equiv $ 1 mod(4) cases read:} \\
\multicolumn{2}{|c|}{} \\ \hline
$C_4$ action on $(C_p \times  C_p)$ &    automorphism group  \\ \hline
$a^p=b^p=(a,b)=c^4=a^c*a^x =(b,c)= \cdots $&$  D_4 \times $Hol$(C_p) \times
C_{(p-1)} $ \\
  $a^p=b^p=(a,b)=c^4=a^c*a^x =b^c*b= \cdots $&$  D_4 \times $ Hol$(C_p) \times $
Hol$(C_p)$     \\
  $ a^p=b^p=(a,b)=c^4=a^c*a^x =b^c*b^y= \cdots  $&$  D_4 \times $[Hol$(C_p)wrC_2]$   \\
    \quad or equivalently \qquad  $=a^c*b^{-1}=b^c*a= \cdots $&
                 \\
   $ a^p=b^p=(a,b)=c^4=a^c*a^x=b^c*b^x= \cdots $&$  D_4 \times $ Hol$(C_p \times
C_p)  $  \\ \hline
\end{tabular}\\
\linebreak\\

\begin{tabular}{|c|c|c|} \hline
\multicolumn{3}{|c|}{Table 4b} \\ \hline
\multicolumn{3}{|c|}{  $p\equiv  1$ mod(4) cases } \\
\multicolumn{3}{|c|}{Case of $C_4 \times  C_2$ actions } \\ \hline
  2-group     &   $  C_4 \times  C_2$ action(s)   &       automorphism \\
     &       on $(C_p \times  C_p)$     &         group(s)  \\ \hline
$ C_8 \times  C_2$ & [$a,b$] [$ab,b$] [$a,ab$]    &$ C_2 \times  C_2 \times
$Hol$(C_p) \times $ Hol$(C_p)$ \\ \hline
    $ <\!2,2|2\!>$     & [$a,b$] [$ab,b$] [$a,ab$]   &$ C_2 \times  C_2 \times$
Hol$(C_p) \times $ Hol$(C_p)$ \\ \hline
    $ C_4 \times  C_4$     &    [$a,b$]       &$ C_2 \times  C_2 \times$
Hol$(C_p) \times $ Hol$(C_p)$ \\ \hline
    $ C_4 \times  C_4$      & [$a,ab$]    &$ C_2 \times  C_2 \times $
(Hol$(C_p)wrC_2$) \\ \hline
  $ C_4 \times  C_2 \times  C_2 $&     [$a,b$]    &$ C_2 \times  C_2 \times$
Hol$(C_p) \times $ Hol$(C_p)$ \\ \hline
    $ C_4 \times  C_2 \times  C_2 $&      [$a,ab$]  &$ C_2 \times  C_2 \times
$ (Hol$(C_p)wrC_2$) \\ \hline
$ C_4 @ C_4 $  &      [$a,b$]  $\ast$  &$ C_2 \times  C_2 \times $ (Hol$(C_p)
\times $Hol$C_p$) \\ \hline
$C_4 @ C_4   $&      [$a,ab$] $\ast$  &$ C_2 \times  C_2 \times $ (Hol$(C_p)wrC_2$) \\
\hline
$(4,4|2,2)$   &      [$a,c$]      &$ C_2 \times  C_2 \times$  Hol$(C_p) \times
$ Hol$(C_p)$ \\ \hline
  $(4,4|2,2)$   &      [$a,ac$] &$ C_2 \times  C_2 \times $ (Hol$(C_p)wrC_2$) \\
\hline
\multicolumn{3}{|c|}{$ \ast$ Here the element $b$ is acting as a $C_4$ on
the $p$-group.         } \\ \hline
\multicolumn{3}{|c|}{Notation: Relations for $(4,4|2,2) [a,ac]$ are:
                } \\
\multicolumn{3}{|c|}{} \\
\multicolumn{3}{|c|}{$a^4=b^2=c^2=(a,b)=(b,c)=a^c*a^{-1}*b=$} \\
\multicolumn{3}{|c|}{$d^p=e^p=(d,e)= $} \\
\multicolumn{3}{|c|}{$d^a*d^x=e^a*e^x=$} \\
\multicolumn{3}{|c|}{ $(d,b)=(e,b)= $} \\
\multicolumn{3}{|c|}{$(d,c)=e^c*e=1$} \\
\multicolumn{3}{|c|}{The rest follow in a like manner. For $C_8 \times  C_2$
and $<\!2,2|2\!>$ all three } \\
\multicolumn{3}{|c|}{extensions for these two 2-groups have the same
automorphism group.} \\ \hline
\end{tabular}

\begin{tabular}{|l|c|c|} \hline
\multicolumn{3}{|c|}{Table 5a} \\
\multicolumn{3}{|c|}{Groups of order $16p^2$ without a normal sylow $p$-subgroup}\\
\multicolumn{3}{|c|}{$p=3$ cases with a normal sylow 2-subgroup }\\ \hline
\qquad \quad group $g$    &  &                             Aut($g$) \\ \hline
1. [$(C_4 \times  C_4) @ C_3] \times  C_3$& &        $  S_3 \times  [384]$ \\
2. $(C_4 \times  C_4) @ C_9 $&   $      C_3$ action & $  C_3 \times $
[384] \\
3. $A_4 \times  C_4 \times  C_3 $ & &      $   S_4 \times  S_3 \times  C_2 $
\\
4. $(1^2 @ C_9) \times  C_4$&$ C_3$ action &$    S_4 \times  C_3 \times  C_2 $ \\
5. $ A_4 \times  A_4$ & &               $    S_4wrC_2$      \\
6. $(1^4 @ C_3)  \times  C_ 3 $ & &            $ S_3 \times  [5760]$ \\
7. $1^4 @ C_9$ &$            C_3$ action &$    C_3 \times  [5760] $  \\
8. $ A_4 \times  C_2 \times  C_2 \times  C_3$ & & $   S_4 \times  S_3 \times
S_3$  \\
9. $(C_2 \times  C_2) @ C_9  \times  C_2 \times  C_2$ & $  C_3$ action &$
S_4 \times  S_3 \times  C_3$  \\
10. $ SL(2,3) \times  C_2 \times  C_3$& & $  S_4 \times  S_3 \times  C_2 $
\\
  11. $(Q_2 @ C_9) \times  C_2$ &$      C_3$ action &$  S_4 \times  C_3
\times  C_2 $ \\
  12. $(C_4$Y$Q_2) @ C_3 \times  C_3$ & &  $    S_4 \times  S_3 \times  C_2
$ \\
  13. $(C_4$Y$Q_2) @ C_9 $&$      C_3$ action &$    S_4 \times  C_3 \times
C_2$  \\ \hline
\multicolumn{3}{|c|}{} \\
\multicolumn{3}{|c|}{$1^2$ means the group  $C_2 \times  C_2. $ } \\
\multicolumn{3}{|c|}{$1^4$ means the elementary abelian group of order 16.}
\\
\multicolumn{3}{|c|}{} \\
\multicolumn{3}{|c|}{The groups [384] and [5760] are the same complete
groups   } \\
\multicolumn{3}{|c|}{as found in Table 2b above. The group $S_4wrC_2$ is
also a  } \\
\multicolumn{3}{|c|}{ complete group of order 1152.
    } \\
\multicolumn{3}{|c|}{} \\ \hline
\end{tabular}
\linebreak\\

\begin{tabular}{|l|c|c|} \hline
\multicolumn{3}{|c|}{Table 5b} \\
\multicolumn{3}{|c|}{Groups of order $16p^2$ without a normal sylow $p$-subgroup } \\
\multicolumn{3}{|c|}{$p=3$ cases without any normal sylow subgroup} \\ \hline
\multicolumn{3}{|c|}{Direct-product cases } \\ \hline
\qquad \quad group $g$   &   &   Aut($g$) \\ \hline
1. $S_4 \times  S_3$ & &  $      S_4 \times  S_3 $        \\
2. $S_4 \times  C_2 \times  C_3 $& &  $ S_4 \times  C_2 \times  C_2 $ \\
3. $S_3 \times  SL(2,3)$ & & $ S_4 \times  C_2 \times  C_2 $     \\
4. $A_4 \times  S_3 \times  C_2 $& & $ S_4 \times  C_2 \times  C_2 $     \\
5. $(C_2 \times  C_2) @ D_9 \times  C_2$& & $    C_2 \times  [216] $ \\
6. $<2,3,4> \times  C_3$ & & $   S_4 \times  C_2 \times  C_2$  \\
7. $GL(2,3) \times  C_3$ &  & $    S_4 \times  C_2 \times  C_2$  \\
8. $(3,3,3;4)) \times  C_2 $ &    & $     C_2 \times  [432]$  \\
9. $ (A_4 @ C_4) \times  C_3$ & $C_4$ acts on $A_4$ as a $C_2$&$  S_4
\times  C_2 \times  C_2$ \\
10. $A_4 \times  Q_3$  & &   $      S_4 \times  S_3 \times  C_2$ \\ \hline
\multicolumn{3}{|c|}{Non-direct-product cases} \\ \hline
\multicolumn{2}{|c|}{group $g$  } &   Aut($g$) \\ \hline
\multicolumn{2}{|l|}{11. $[(C_2 \times  C_2) @ C_9] @ C_4 $: }& $   C_2
\times  [216]$   \\
\multicolumn{2}{|c|}{$ a^4=b^2=c^2=p^9=a^b*c*a^{-1}= $} & \\
\multicolumn{2}{|c|}{$(a,c)=(b,c)=b^p*c=c^p*c*b=p^a*p=1$; } & \\
\multicolumn{2}{|l|}{12. } &$C_2 \times [216] $\\
\multicolumn{2}{|c|}{$a^8=b^2=p^9=a^b*a^5=(a^2)^p*a^{-1}*b=$ }& \\
\multicolumn{2}{|c|}{$(b*a)^p*a*b=a^{-1}*p*a^6*b*p=1$;     }& \\
\multicolumn{2}{|l|}{13. } &$C_2 \times [216] $\\
\multicolumn{2}{|c|}{$ a^8=b^4=p^9=a^4*b^2=a^b*a= $}& \\
\multicolumn{2}{|c|}{$(a^2)^p*b^3=b^p*b^3*a^6=       $ }& \\
\multicolumn{2}{|c|}{$a^{-1}*p*a^3*p=1$;      }& \\
\multicolumn{2}{|l|}{14. $[A_4 \times  C_3] @ C_4$: }& $  C_2 \times  [432]
$ \\
\multicolumn{2}{|c|}{$ a^4=b^2=c^2=p^3=q^3=a^b*c*a^{-1} = $ }&\\
\multicolumn{2}{|c|}{$(a,c)=(b,c)=(p,q)=q^a*q=b^q*c= $ } & \\
\multicolumn{2}{|c|}{$c^q*c*b=p^a*p=(b,p)=(c,p)=1$; } & \\
\multicolumn{2}{|l|}{15.} &$C_2 \times [432] $\\
\multicolumn{2}{|c|}{$a^8=b^4=a^4*b^2=a^b*a=p^3=q^3=(p,q)= $}& \\
\multicolumn{2}{|c|}{$(a^2)q*b^{-1}=b^q*b^{-1}*a^6=$ }& \\
\multicolumn{2}{|c|}{$ a^{-1}*q*a^3*q=p^a*p=(b,p)=1$;  }& \\
\multicolumn{2}{|l|}{16. } &$C_2 \times [432] $\\
\multicolumn{2}{|c|}{$a^8=b^2=a^b*a^5=p^3=q^3=(p,q)=$ }& \\
\multicolumn{2}{|c|}{$(a^2)^q*a^{-1}*b=(b*a)^q*a*b=$  }& \\
\multicolumn{2}{|c|}{$ a^{-1}*q*a^6*b*q=p^a*p=p^b*p=1$; }& \\
\multicolumn{2}{|l|}{17.} &$S_4 \times S_3 \times C_2 $\\
\multicolumn{2}{|c|}{$a^4=a^2*b^2=a^b*a=c^3=a^c*b^{-1}= $ }& \\
\multicolumn{2}{|c|}{$ b^c*b^{-1}*a^{-1}=d^4=f^3=f^d*f=$  }& \\
\multicolumn{2}{|c|}{$ a^2*d^2=(a,f)=(b,f)=(c,f)=(b,d)=(c,d)=1$.}& \\
\multicolumn{2}{|c|}{This group comes from identifying the } &\\
\multicolumn{2}{|c|}{centers of the two groups $Q_2 @ C_3$  and $ C_3 @ C_4$.}& \\ \hline
\end{tabular} \\
\linebreak\\

\begin{tabular}{|c|} \hline
Notes for Table 5b \\ \hline
            A presentation for the above occurring (complete) group of \\
                     order 216 is:                                     \\
                                                                       \\
             $a^2=b^2=c^3=(a,c)=(a*d)^2=(b,c)=c*d^2*c^{-1}*d=$           \\
             $a*b*a*d*b*d^{-1}=(a*b)^2*d^{-1}*b*d=1$.                    \\
                                                                       \\
               The above group of order 432 is a complete group with      \\
               20 conjugacy classes and is not Hol$(C_3 \times  C_3)$.  \\
               A presentation for this group is:                       \\
                                                                       \\
            $ a^4=b^2=c^3=d^2=(a,c)=(a,d)=(a,d)=(b,d)=(c*d)^2=(a*b)^3 $\\
                     $=(b*c)^2*(b*c^{-1})^2=a*b*a*c*b*c*a*b*c=1. $ \\ \hline

\end{tabular}
\linebreak\\

\begin{tabular}{|l|c|c|} \hline
\multicolumn{3}{|c|}{Table 5c} \\ \hline
\multicolumn{3}{|c|}{} \\
\multicolumn{3}{|c|}{Groups of order $16p^2$ without a normal sylow
$p$-subgroup     } \\
\multicolumn{3}{|c|}{} \\
\multicolumn{3}{|c|}{Cases with $p > 3$} \\
\multicolumn{3}{|c|}{} \\ \hline
\multicolumn{3}{|c|}{ $p=5$: 2 cases } \\ \hline
1. $(1^4 @ C_5) \times  C_5 $   & &    [960]$ \times $ Hol$(C_5)$ \\
2. $1^4 @ C_{25}$ & $ C_5$ action &      [960]$ \times  C_5 $ \\ \hline
\multicolumn{3}{|c|}{See Table 2c for relations of the group [960]
(complete). } \\ \hline
\multicolumn{3}{|c|}{ $p=7$: 3 cases } \\ \hline
1. $(1^3 @ C_7) \times  C_2 \times  C_7 $ & &       [168]$ \times$
Hol$(C_7)$ \\
2. $(1^3 @ C_{49}) \times  C_2$  & $ C_7$ action &    [168]$ \times  C_7 $
\\
3. nonnormal case:      & & \\
$(1^3 @ C_7) \times  D_7$ & &          [168]$ \times $ Hol$(C_7)$ \\ \hline
\multicolumn{3}{|c|}{See Table 2c for relations of the group [168]
(complete). } \\ \hline
\end{tabular}
\linebreak\\

\begin{tabular}{|c|c|} \hline
\multicolumn{2}{|c|}{Table 6a} \\ \hline
\multicolumn{2}{|c|}{} \\
\multicolumn{2}{|c|}{Number of groups of order $16p^2$} \\
\multicolumn{2}{|c|}{(Taken from  R. Nyhlen \cite{5})} \\
\multicolumn{2}{|c|}{} \\ \hline
     $p= 3$      &        197                      \\
                   $p= 5$              &221                      \\
                   $p= 7$              &172                       \\
            $\ast$   $p\equiv  1$ mod(16)-----$>$&257 [e.g., $p=17$]        \\
             \quad   $p\equiv 3$ mod(8)------$>$&167 [e.g., $p=11$, 19]      \\
             \quad      $p\equiv 5$ mod(8)------$>$&219 [e.g., $p=13$, 29, 37]  \\
              \quad    $p\equiv  7$ mod(8) -----$>$&169 [e.g., $p=23$, 31]     \\
             \quad     $p\equiv  9$ mod(16)-----$>$&243 [e.g., $p=41$]  \\ \hline
\multicolumn{2}{|c|}{} \\
\multicolumn{2}{|c|}{$\ast$ Number here might be 258; see Table 6b for
breakdown.} \\\multicolumn{2}{|c|}{} \\
\multicolumn{2}{|c|}{There appear to be 15 different groups arising from }\\
\multicolumn{2}{|c|}{
             the extension $ (C_{17} \times  C_{17}) @ C_{16}$. In his thesis
Nyhlen }\\
\multicolumn{2}{|c|}{ states there are 14 cases. R. Laue gets 15 cases } \\
\multicolumn{2}{|c|}{(private communication, Feb. 15, 1994). }\\ \hline
\end{tabular}
\linebreak\\

\begin{tabular}{|c|c|c|c|c|c|} \hline
\multicolumn{6}{|c|}{Table 6b} \\
\multicolumn{6}{|c|}{Additional groups of order $16p^2$ for $p> 3$ }\\
\hline
2-group & image of & \multicolumn{4}{|c|}{prime $p$} \\ \cline{3-6}
  & 2-group& 5 mod(8) & 7 mod(8) & 9 mod(16) & 1 mod(16) \\ \hline
      $ C_{16}$  &$   C_4$    &   3 + 1  &        &         &    \\
                &$   C_8$      &           &    &  7 + 1   &
  \\
                &$   C_{16}$    &        &  1     &     1 &       14 +
1    \\ \hline
  $  C_8 \times  C_2 $ &$   C_4 $   &   6 + 2 &     none &       &
none  \\
             &$ C_4 \times  C_2$   &     3  &     &               &
   \\
        &$   C_8$     &    &    &          7 + 1 &               \\
      & $C_8 \times  C_2$ &  &    &            4 &            \\
\hline
   $ C_4 \times  C_4 $ &$    C_4$   &   3 + 1  &    none &      none &
none  \\
       & $ C_4 \times  C_2$  &     2  & & &\\
            &$  C_4 \times  C_4$  &   1 & & &
\\ \hline
$C_4 \times  C_2 \times  C_2$  &$   C_4 $  &   3 + 1  &    none  &   none  &
none  \\
          &$    C_4 \times  C_2$  &     2 & & &                       \\
\hline
   $ C_4 \text{Y} Q_2$ &$  C_4 \text{Y} Q_2$  &     1 &   none &   none  &    none    \\
\hline
$(4,4|2,2)$ &$    C_4 $ &   3 + 1  &    none   &    none &     none   \\
     &$  C_4 \times  C_2  $&     2&  & & \\ \hline
  $ C_4 @ C_4 $ &$    C_4 $    &   3 + 1 &     none &      none &     none
  \\
         &$        C_4 \times  C_2$ &     2 & & & \\ \hline
  $<2,2|2>$ &$    C_4$    &   6 + 2  &    none   &    none &     none
\\
           &$     C_4 \times  C_2$  &    3 & & & \\
                & $<2,2|2>$   &     1 & & & \\ \hline
   $     D_8$     & $    D_8 $   &   none    &    1    &       1  &    none
\\ \hline
     $<-2,4|2>$  &  $<-2,4|2>$ &   none    &   none   &      1  &
none  \\ \hline
   $    Q_4 $   &$    Q_4   $&   none    &    1    &       1   &    none
\\ \hline
        totals  &          &    53   &      3    &      24   &     15   \\
\hline
\multicolumn{6}{|c|}{ Notation: The numbers 3 + 1 (e.g., in the $C_4$ images
for $C_{16}$) mean  } \\
\multicolumn{6}{|c|}{that we get 3 additional groups with a $C_4$ action
when } \\
\multicolumn{6}{|c|}{$p \equiv 5 $ mod(8) for extensions of the form            } \\
\multicolumn{6}{|c|}{} \\
\multicolumn{6}{|c|}{$(C_p \times  C_p) @ C_{16}$                 } \\
\multicolumn{6}{|c|}{} \\
\multicolumn{6}{|l|}{and one case from the extension
$(C_{p^2})@ C_{16}$, and similarly for the      } \\
\multicolumn{6}{|l|}{other cases.} \\
\multicolumn{6}{|c|}{The $p\equiv  5$ mod(8), 9 mod(16) and 1 mod(16) cases here
        } \\
\multicolumn{6}{|l|}{represent the same type of situation found in the $16p$
cases for which } \\
\multicolumn{6}{|l|}{$p\equiv 1$ mod(4) only, $p\equiv 1$ mod(8) but not 1 mod(16), and
$p \equiv  1$ mod(16)  } \\
\multicolumn{6}{|l|}{respectively. None means no additional groups for this
prime, but if  } \\
\multicolumn{6}{|l|}{$p \equiv  9$ mod(16) (or  1 mod(8)), then the groups for $p\equiv  5$
mod(8) also    } \\
\multicolumn{6}{|l|}{have recurrences in the primes $p\equiv  9$ mod(16) and $p\equiv  1$
mod(16), and   } \\
\multicolumn{6}{|l|}{likewise for the $p\equiv  9$ mod(16) and $p\equiv  1$ mod(16) cases.
               } \\ \hline
\end{tabular}

\begin{tabular}{|c|c|c|} \hline
\multicolumn{3}{|c|}{Table 7} \\
\multicolumn{3}{|c|}{} \\
\multicolumn{3}{|c|}{Groups and automorphism groups for $(C_p \times  C_p) @
(C_8 \times  C_2)$ } \\ \hline
$a^p=b^p=(a,b)=c^8=d^2=(c,d)= $ & & automorphism group (?) \\ \hline
\multicolumn{3}{|c|}{} \\ \hline
$a^c*a^x=(b,c)=(a,d)=b^d*b=1;$&$  [C_p @ C_8 \times  D_p]$ &   Hol$(C_p)
\times $ Hol$(C_p)  \ast$ \\ \hline
$a^c*a^x=b^c*b^t=(a,d)=b^d*b=1;$&$  t = x^2$ &   Hol$(C_p) \times $
Hol$(C_p)  \ast$ \\ \hline
$a^c*a^x=b^c*b^t=a^d*a=(b,d)=1;$&$  t = x^3$ &   Hol$(C_p)wrC_2  \ast$
\\ \hline
$a^c*a^x=b^c*b^t=a^d*a=(b,d)=1;$&$  t = x^5$ &   Hol$(C_p)wrC_2  \ast$
\\ \hline
\multicolumn{2}{|c|}{where $x^8 \equiv  1$ mod($p$)}& $\ast$ verified for $p= 17$
\\ \hline
\end{tabular} \\
\linebreak\\

\begin{tabular}{|c|c|c|c|c|c|c|c|c|} \hline
\multicolumn{9}{|c|}{Table 8} \\
\multicolumn{9}{|c|}{} \\
\multicolumn{9}{|c|}{Groups and automorphism groups for $(C_{17} \times
C_{17}) @ C_{16}$ } \\ \hline
   $x =$   &   1  &   2 &     3 &    4  &    5 &    6 &    7  &   8  \\ \hline
ncl($g$) =  & 97 &  37  &     34  &   49 &   34 &   34 &   34  &  37 \\ \hline
Aut($g$) = &(8,2) & (8,2) & Hol$(C_p \times  C_p)$& (8,2) & (8,2) &(9,2)&
(8,2)& (8,2) \\ \hline
\multicolumn{9}{|c|}{} \\ \hline
$x=$ &  9  &  10 &    11 &    12  &   13   &  14  &   15  &   16 \\ \hline
ncl($g$) =&   37 &    34 &    34 &    34 &    49 &    34 &    37 & 289 \\ \hline
Aut($g$) =&  (8,2) & (8,2)&  (9,2) & (8,2) & (8,2) & (9,2) & (8,2) & (8,1) \\
\hline
\multicolumn{9}{|c|}{Here ($a,b$) stands for $2^a*17^b$.   All of the groups
in the Aut($g$)  } \\
\multicolumn{9}{|c|}{ row labelled (8,2) are the groups Hol$(C_{17}) \times
$ Hol$(C_{17})$. The ones } \\
\multicolumn{9}{|c|}{labelled (9,2) are [Hol$(C_{17})wrC_2$], and the
last one (8,1) is}\\
\multicolumn{9}{|c|}{Hol$(C_{17}) \times  C_{16}$. ncl($g$) is the number of
conjugacy classes  } \\
\multicolumn{9}{|c|}{in the group $g$. }\\ \hline
\end{tabular}\\

\begin{tabular}{|c|c|c|c|} \hline
\multicolumn{4}{|c|}{Table 9} \\
\multicolumn{4}{|c|}{New groups of order $16p^2$ with $p\equiv  7$ mod(8) or 9
mod(16)}\\
\multicolumn{4}{|c|}{}\\ \hline
2-group & type & automorphism group & comments \\ \hline
$C_{16} $&   7  & $ (C_7 \times  C_7) @ (C_3 \times  QD_8)$ &     \\
        &  7  &$ (C_{23} \times  C_{23}) @ (C_{11} \times  (C_3 @ QD_8))$&
$\ast$ \\
        &  7  & $(C_{31} \times  C_{31}) @ (C_3 \times  C_5 \times
QD_{64})$& \\ \hline
$D_8$ and $Q_4$ &  7 &$ (C_7 \times  C_7) @ ( C_3 \times  QD_{16} )$&  \\
                 &  7 &$  (C_{23} \times  C_{23}) @ (C_{11} \times  QD_{16}
)$&   \\
  $D_8$   &  7  &$  (C_{31} \times  C_{31}) @ (C_3 \times  C_5 \times
D_{16})$& \\
  $Q_4$   &  7  &$  (C_{31} \times  C_{31}) @ (C_3 \times  C_5 \times
Q_8)$& \\ \hline
$ C_{16}$ &  9  &   [ Hol$(C_p)wrC_2$ ] conjecture  & $\ast\ast$  not
done  \\ \hline
  $D_8$ and $Q_4$ & 1,9 &  $  (C_{17} \times  C_{17}) @ (C_{16} \text{Y} D_{16} )
$&$ \ast\ast\ast$   \\
  $QD_8$   & 1,9 &$ (C_{17} \times  C_{17}) @ (C_{16} \text{Y} (C_4wrC_2))$  &
\\ \hline
\multicolumn{4}{|c|}{$\ast$ The group $(C_3 @ QD_8)$ is isomorphic to
$D_{48}. $}\\
\multicolumn{4}{|c|}{ Is there a higher-order case (e.g., $p= 47$) with a }\\
\multicolumn{4}{|c|}{non-direct-product factor here such as $C_{23} \times
(C_3 @ QD_{32})$ ? }\\ \hline
\multicolumn{4}{|c|}{$\ast\ast$ First case appears for the prime $p= 41$, has
an order that was  }\\
\multicolumn{4}{|c|}{too large to do at this site.}\\ \hline
\multicolumn{4}{|c|}{$\ast\ast\ast$The groups of this form  occur for $p\equiv  1$
or 9 mod(16).    }\\
\multicolumn{4}{|c|}{ The calculations given here are for the case of $p=
17$.   }\\
\multicolumn{4}{|c|}{The next case occurs for $p= 41$, which was too large for
   }\\
\multicolumn{4}{|c|}{us to handle at this site.
}\\
\multicolumn{4}{|c|}{}\\
\multicolumn{4}{|c|}{The groups $C_{16} \text{Y} D_{16}$ and $C_{16} \text{Y} Q_8$ are
isomorphic; so the  }\\
\multicolumn{4}{|c|}{cases $p\equiv  +1$ mod(16) may be expressible as
}\\
\multicolumn{4}{|c|}{}\\
\multicolumn{4}{|c|}{$(C_p \times  C_p) @ ( C_{p-1} \text{Y} D_{16})$ for the $D_8$
and $Q_4$ cases.  }\\
\multicolumn{4}{|c|}{ For $p\equiv  1$ mod(8) (e.g., $p= 41$) the orders of the  }\\
\multicolumn{4}{|c|}{groups involved were too large to handle here.
}\\ \hline
\end{tabular}\\
\linebreak\\

\begin{tabular}{|c|c|c|c|c|c|c|} \hline
\multicolumn{7}{|c|}{Table 10} \\
\multicolumn{7}{|c|}{Class Structure of $(C_3 \times  C_3) @ C_{16} $
[$C_8$ action]}\\
\multicolumn{7}{|c|}{}\\ \hline
order of element & 2 & 3 & 4 & 6 & 8 & 12 \\ \hline
number of elements &  19 &  8  & 132 &  8  &  72 &  48 \\ \hline
number of classes  &  3  &  1  &  6  &  1  &  4  &  2 \\ \hline
\end{tabular}

\begin{tabular}{|c|c|c|c|c|} \hline
\multicolumn{5}{|c|}{Table 11a} \\ \hline
\multicolumn{5}{|c|}{ Matrix Representations for the Presentation(s) $QD_t$
}\\
\multicolumn{5}{|c|}{( $t = 2^{n+1}$ )    } \\ \hline
\multicolumn{5}{|c|}{ $p\equiv  - 1$ mod($2^n$):    $ n > 1 $    }\\
\multicolumn{5}{|c|}{} \\
\multicolumn{5}{|c|}{$a^t = b^2 = a^b*a^{t1} = 1 $  ( $t1 = 2n -1$ )  } \\
\multicolumn{5}{|c|}{} \\
\multicolumn{5}{|c|}{with the matrix representations: } \\
\multicolumn{5}{|c|}{}\\
\multicolumn{5}{|c|}{
$ a=\begin{pmatrix} 0 & 1 \\ y & x \end{pmatrix} \qquad b = \begin{pmatrix} 1 & 0 \\ x & -1\end{pmatrix} $ } \\
\multicolumn{5}{|c|}{} \\
\multicolumn{5}{|c|}{For the cases of $-1$ mod(8) (or $p\equiv  7$ mod(8)) we have }\\
\multicolumn{5}{|c|}{$y = 1$ and  $(x^2 + 2)^2 \equiv  2$ mod($p$).     }\\ \hline
\multicolumn{5}{|c|}{Values of $n$ and $x$ for selected low-order primes are
given below. }\\
\multicolumn{5}{|c|}{ }\\
\multicolumn{5}{|c|}{In every case if $x$ is an allowed entry in the above
matrix, then so  }\\
\multicolumn{5}{|c|}{ is the entry ($p-x$), so for $p= 47$ and 79 below there
are 8 cases.    }\\  \hline
  $n =$  &            4      &     5     &        6     &         7 \\ \hline
group order  &  64 & 128 & 256 & 512 \\ \hline
  $p=$       &      47     &     31     &      191     &     127  \\ \hline
      $x =$    & (1,4,11,18) & (3,5,7,\ldots ) & (5,12,18,\ldots ) &(1,4,5,6,\ldots ) \\
        $y = 1$   &             &            &              &             \\
                 &             &  16 cases  &   32 cases   &   64 cases  \\
\hline
$p =$       &      79     &     223    &              &  \\ \hline
     $x =$     &(8,13,17,24) &(23,29,35,\ldots )&              &             \\
       $y = 1$   &            &            &              &             \\
                 &   8 cases   &  16 cases  &              &             \\
\hline

\multicolumn{5}{|c|}{An algebraic characterization of the parameter $x$ in the
matrix   }\\
\multicolumn{5}{|c|}{ $a$ above is given below, along with an alternate
presentation }\\
\multicolumn{5}{|c|}{for the sylow 2-subgroups of $GL(2,p)$ when $p\equiv -1$
mod(4),\ldots .  }\\ \hline
\end{tabular}

\begin{tabular}{|c|} \hline
Table 11a (continued) \\
\\
Alternate Matrix Representations and Presentation(s) for $QD_t$ \\
($ t = 2^{(n+1)}$) \\ \hline
           An alternate matrix representation for the group $QD_{2^{(n+1)}}$
  is given \\
              by Carter and Fong \cite{8}, namely: \\
$ a = \begin{pmatrix} 0 & 1 \\ 1 & x \end{pmatrix} \qquad b = \begin{pmatrix} 0 & 1 \\ -1 & 0 \end{pmatrix} $
\\
                                                                        \\
           The general expression for this $x$ is \cite{8}:                     \\
                                                                        \\
                               $ x = t + t^p$
\\
                                                                        \\
            where $t$ is a $2^{(n + 1)}$-st root of unity in $GF(p^2)$.
   \\
       A presentation for this alternate representation is:             \\

                                                                        \\
           $  (a*b)^2=(a*(b^{-1}))^2=b^4=a^y*b*a^{-y}*b=1; $   $y= |G|/8.$  \\
  \hline
    \\
                                                                          \\
        The group $<a,b>$ is the sylow 2-subgroup of $GL(2,p)$ with           \\
            $p\equiv  -1$ mod($2^n$) and the $x$ in $<a>$ is a root of the following
  \\
            polynomial (with the appropriate value of $n$):              \\
\multicolumn{1}{|l|}{for }  \\
\multicolumn{1}{|l|}{       $ n = 3 \quad    (x^2 + 2)^2                              \equiv  2$ mod($p$)  } \\
 \multicolumn{1}{|l|}{      $ n = 4 \quad   ((x^2 + 2)^2 - 2)^2                   \equiv  2$ mod($p$)
}\\
 \multicolumn{1}{|l|}{      $ n = 5 \quad  (((x^2 + 2)^2 - 2)^2 - 2)^2                     \equiv  2$ mod($p$)
 } \\
 \multicolumn{1}{|l|}{      $ n = 6 \quad ((((x^2 + 2)^2 - 2)^2 - 2)^2 - 2)^2               \equiv  2$ mod($p$) } \\
\multicolumn{1}{|l|}{       $ n = 7 \quad (((((x^2 + 2)^2 - 2)^2 - 2)^2 - 2)^2 - 2)^2 \equiv  2$
mod($p$) }\\
\multicolumn{1}{|l|}{etc. }  \\
\hline
        Matrix representations for the sylow 2-subgroups of $GL(n,p)$ are   \\
                            discussed in Appendix 3.                \\
\hline
\end{tabular}

\begin{tabular}{|c|} \hline
Table 11b \\ \hline
Matrix Representations for the Presentation(s) $C_twrC_2$ $(t=2^n)$:\\
\\
$a^t = b^2 = (a*b)^2=(a^{-1}*b)^2=1;$\\
\\
for the primes $p\equiv  +1$ mod($2^n$), $n > 1.$       \\
     Some choices here are: \\
 $ a = \begin{pmatrix} 1 & 0 \\ 0 & x \end{pmatrix} \quad b=\begin{pmatrix} 0 & 1 \\ 1 & 0\end{pmatrix}, $
   \\
                                                          \\
        where here $x$ is a $2^n$-th root of unity in $GF(p)$.                 \\
\hline

           Note that if $n = 4$, then this group is $C_{16}wrC_2$ .  In
many \\
       cases we would like a representation for $C_8wrC_2$ or $C_4wrC_2$ over \\
       the field $GF(17)$. In these cases we have:
          \\
               $   h =$ $<a^2,b>$ is just $C_8wrC_2$,   \\
            \\
                  $k =$ $<a^4,b>$ is just $C_4wrC_2$.    \\
                                                    \\
    Hence knowing a matrix representation for $C_twrC_2 (t=2^n)$ will also
   \\
     yield matrix representations for other wreath products of interest. \\
\hline

\end{tabular}

\begin{thebibliography}{99}

      \bibitem{1} Becker, W., \textquotedblleft A Preliminary Report of a Computer Study of
            Finite Groups and their Automorphism Groups", unpublished
            manuscript (1994). This report is rather long (200 pages)
            and contains results on groups of orders $p^3q$, $p^4q$,
            $p^3q^{2}$ (incomplete results here), orders 240, $32p$, $36p$
            ($p > 3$), and groups of order $64p$ without a normal
            odd order sylow subgroup, as well as a long list of complete
            groups which were found in this study. The groups of order
            $8p^2$ were not included in this report but are reported in \cite{2}.\\

\bibitem{2} Becker, W. and Becker, Elaine W., \textquotedblleft The automorphism groups of the groups of order $8p^2$", arXiv.math.GR0610555.\\

\bibitem{3} Coxeter, H.S.M. and Moser, W.O.J., \textquotedblleft Generators and Relations for Discrete
Groups", second edition, Springer-Verlag, New York, 1965.\\

\bibitem{4} Sag, T. W. and Wamsley, J. W., \textquotedblleft Minimal Presentations for Groups of Order $2^n $, $n \leq  6$". J. of the Australian Math. Soc., vol. 15, 1973, pp. 461--469.\\

\bibitem{5} Nyhlen, R., \textquotedblleft Determination of the abstract groups of order $16p^2$ and $8p^3$", Unpublished thesis, Uppsala, Sweden, 1919.\\

\bibitem{6} Becker, W. \textquotedblleft Automorphism groups for the groups of the form ($C_{p} \times  C_{p}$) @ $X$, where $X$ is a 2-group and has action $D_4$ or $Q_{2}$ on the $p$-group", in preparation.\\

\bibitem{7} Senior, J. K. and Hall, M., \textquotedblleft The Groups of Order $2^n$, $n \leq
6$", Macmillan, New York, 1964.\\

\bibitem{8} Carter, R. and Fong, P. \textquotedblleft The Sylow 2-Subgroups of the Finite Classical Groups" J. of Algebra, vol. 1, 1964, pp. 139--151 (see p. 143).\\

\bibitem{9} This appendix is taken from the appendix in \cite{2}.\\

\bibitem{10} Relations of the form $x^s \equiv  1$ mod($p$) (or $z$, or $t$,\ldots) occur in the text and the appendices. In each case $x$ is an $s^{th}$ root of unity. (Here $s = 3$, 4,\ldots , $q$,\ldots,
$p-1$.)\\

\bibitem{11} Lunn, A. C. and Senior, J. K., \textquotedblleft A Method of Determining all the
Solvable Groups of Given Order and its Application to the Orders 16$p$ and 32$p$", Amer. J. Math.,
vol. 56, 1934, pp. 319--327.\\

\bibitem{12} Becker, W. and Becker, Elaine W., \textquotedblleft On the groups and automorphism groups of the groups of order 64$p$ without a normal sylow $p$-subgroup", to appear.\\

\bibitem{13} Besche, Hans Ulrich, Eick, Bettina, and O'Brien, E. A., \textquotedblleft A millennium
project: Constructing small groups", International Journal of Algebra and Computation, vol. 12,
no. 5, 2002, pp. 623--644. This is just one of several articles dealing with the constructing
and the listing of the results of their calculations.


\end{thebibliography}
\end{document}